\documentclass{article}

\usepackage[left=32mm,right=32mm,top=32mm,bottom=32mm]{geometry}
\usepackage{color}
\usepackage{amsthm,amsmath,amssymb}
\usepackage{booktabs}
\usepackage{mathpazo}
\usepackage{microtype}
\usepackage{overpic}
\usepackage{bm}
\usepackage{sectsty}
\usepackage[
	pdftitle={Moduli spaces of hyperbolic surfaces and their Weil-Petersson volumes},
	pdfauthor={Norman Do},
	ocgcolorlinks,
	linkcolor=linkblue,
	citecolor=linkred,
	urlcolor=linkred]
{hyperref}
\usepackage[hang,flushmargin]{footmisc}
\usepackage{enumitem}
\usepackage{titlesec}
	\titlespacing{\section}{0pt}{12pt}{0pt}
	\titlespacing{\subsection}{0pt}{6pt}{0pt}

\makeatletter 

\long\def\@footnotetext#1{%
\H@@footnotetext{%
\ifHy@nesting 
\hyper@@anchor{\@currentHref}{#1}%
\else 
\Hy@raisedlink{\hyper@@anchor{\@currentHref}{\relax}}#1%
\fi 
}}

\def\@footnotemark{%
\leavevmode 
\ifhmode\edef\@x@sf{\the\spacefactor}\nobreak\fi 
\H@refstepcounter{Hfootnote}%
\hyper@makecurrent{Hfootnote}%
\hyper@linkstart{link}{\@currentHref}%
\@makefnmark 
\hyper@linkend 
\ifhmode\spacefactor\@x@sf\fi 
\relax 
}%

\ifFN@multiplefootnote%
\renewcommand*\@footnotemark{%
\leavevmode 
\ifhmode 
\edef\@x@sf{\the\spacefactor}%
\FN@mf@check 
\nobreak 
\fi 
\H@refstepcounter{Hfootnote}%
\hyper@makecurrent{Hfootnote}%
\hyper@linkstart{link}{\@currentHref}%
\@makefnmark 
\hyper@linkend 
\ifFN@pp@towrite 
\FN@pp@writetemp 
\FN@pp@towritefalse 
\fi 
\FN@mf@prepare 
\ifhmode\spacefactor\@x@sf\fi 
\relax%
}%
\fi 

\makeatother 

\theoremstyle{plain}
\newtheorem{theorem}{Theorem}
\newtheorem{proposition}[theorem]{Proposition}
\newtheorem{corollary}[theorem]{Corollary}
\newtheorem{lemma}[theorem]{Lemma}
\newtheorem{conjecture}[theorem]{Conjecture}

\theoremstyle{definition}

\newtheorem{example}[theorem]{Example}

\definecolor{linkred}{rgb}{0.6,0,0}
\definecolor{linkblue}{rgb}{0,0,0.6}

\newcommand{\LL}{{\mathbf L}}
\newcommand{\M}{\overline{\mathcal M}}
\newcommand{\T}{{\mathcal T}}
\newcommand{\x}{\mathbf{x}}

\sectionfont{\large}
\subsectionfont{\normalsize}

\begin{document}

{\large \bfseries Moduli spaces of hyperbolic surfaces and their Weil--Petersson volumes}

{\bfseries Norman Do}

{\em Abstract.} Moduli spaces of hyperbolic surfaces may be endowed with a symplectic structure via the Weil--Petersson form. Mirzakhani proved that Weil--Petersson volumes exhibit polynomial behaviour and that their coefficients store intersection numbers on moduli spaces of curves. In this survey article, we discuss these results as well as some consequences and applications.

\vspace{-26pt}
\setlength{\parskip}{0pt}
\makeatletter \renewcommand{\@dotsep}{10000} \makeatother
\renewcommand\contentsname{}
\tableofcontents
\setlength{\parskip}{6pt}

\vfill

\rule{0.25\textwidth}{0.5pt} \\
Department of Mathematics and Statistics, The University of Melbourne, Victoria 3010, Australia \\
{\em Email:} \href{mailto:normdo@gmail.com}{normdo@gmail.com}


\newpage

\section{Introduction} \label{introduction}

Over the last few decades, moduli spaces of curves have become increasingly important objects of study in mathematics. In fact, they now lie at the centre of a rich confluence of seemingly disparate areas such as geometry, topology, combinatorics, integrable systems, matrix models and theoretical physics.

Smooth curves with marked points are equivalent to Riemann surfaces with punctures, and these are in turn equivalent to punctured surfaces with complete constant curvature Riemannian metrics.\footnote{We decree that all surfaces referred to in this article are to be connected and oriented. We also decree that all algebraic curves referred to in this article are to be complex, connected and complete.} For all but finitely many pairs $(g,n)$, a genus $g$ surface with $n$ punctures admits a hyperbolic metric. Thus, we are led to the study of hyperbolic surfaces and their corresponding moduli spaces. In this article, we adopt such a hyperbolic geometric perspective, which allows for lines of thought that have no natural analogue in the realms of algebraic geometry and complex analysis.

For an $n$-tuple $\LL = (L_1, L_2, \ldots, L_n)$ of positive real numbers, let ${\mathcal M}_{g,n}(\LL)$ denote the set of genus $g$ hyperbolic surfaces with $n$ geodesic boundary components whose lengths are prescribed by $\LL$. We require the boundary components to be labelled from 1 up to $n$ and consider hyperbolic surfaces up to isometries which preserve these labels. The Teichm\"{u}ller theory construction of this space endows it with an orbifold structure and local coordinates known as Fenchel--Nielsen coordinates. These may be used to define the Weil--Petersson symplectic form $\omega$, thus providing the moduli space of hyperbolic surfaces ${\mathcal M}_{g,n}(\LL)$ with a natural symplectic structure. Our primary focus will be on the Weil--Petersson volume
\[
V_{g,n}(\LL) = \int_{{\mathcal M}_{g,n}(\LL)} \frac{\omega^{3g-3+n}}{(3g-3+n)!}.
\]

The foundation of this article is a result due to Mirzakhani~\cite{mir1} which states that the Weil--Petersson volume is given by the following polynomial.
\[
V_{g,n}(\LL) = \sum_{|\bm{\alpha}| + m = 3g-3+n} \frac{(2\pi^2)^m \int_{\M_{g,n}} \psi_1^{\alpha_1} \psi_2^{\alpha_2} \cdots \psi_n^{\alpha_n} \kappa_1^{m}}{2^{|\bm{\alpha}|} \alpha_1! \alpha_2! \cdots \alpha_n! m!} L_1^{2\alpha_1} L_2^{2\alpha_2} \cdots L_n^{2\alpha_n}
\]
Here, $\psi_1, \psi_2, \ldots, \psi_n \in H^2(\M_{g,n}; \mathbb{Q})$ are the psi-classes while $\kappa_1 \in H^2(\M_{g,n}; \mathbb{Q})$ is the first Mumford--Morita--Miller class on the Deligne--Mumford compactification of the moduli space of curves. We use $\bm{\alpha} = (\alpha_1, \alpha_2, \ldots, \alpha_n)$ to denote an $n$-tuple of non-negative integers and $|\bm{\alpha}|$ to denote the sum $\alpha_1 + \alpha_2 + \cdots + \alpha_n$. For a concise exposition of intersection theory on moduli spaces of curves, including the definitions of the psi-classes and the Mumford--Morita--Miller classes, see Appendix~\ref{intersection-theory}.

Mirzakhani originally demonstrated the polynomiality of Weil--Petersson volumes using symplectic reduction~\cite{mir1}. She subsequently provided an alternative proof where the main idea is to unfold the integral on ${\mathcal M}_{g,n}(\LL)$ to a more tractable cover over the moduli space~\cite{mir2}. One of the tools required is a generalisation of McShane's identity which, in its original form, states that a certain sum over the simple closed geodesics on a hyperbolic once-punctured torus is constant. The end result is a recursive formula which can be used to calculate all Weil--Petersson volumes. This in turn yields a recursion for intersection numbers on $\M_{g,n}$. As a consequence, Mirzakhani was able to deduce the celebrated Witten--Kontsevich theorem.

Some of the recent work on moduli spaces of hyperbolic surfaces has focused on the behaviour of Weil--Petersson volumes under various limits. For example, the Weil--Petersson volume polynomials exhibit interesting behaviour when one of the lengths formally approaches $2\pi i$. This manifests as non-trivial relations between $V_{g,n+1}(\LL, L_{n+1})$ and $V_{g,n}(\LL)$, which generalise the string and dilaton equations~\cite{do-nor}. These can be proven using algebro-geometric arguments, but also indicate that a proof may entail the geometry of hyperbolic cone surfaces. In particular, one can interpret the evaluation $L = \theta i$ as the degeneration of the corresponding geodesic boundary component to a cone point with angle $\theta$.

A direct implementation of Mirzakhani's recursive formula produces a rather slow algorithm for the calculation of Weil--Petersson volumes. However, Zograf has managed to develop an empirically much faster algorithm and calculated enough numerical data to produce interesting conjectures concerning Weil--Petersson volumes in the large $g$ limit~\cite{zog1}. Some progress on these conjectures has recently been made by Mirzakhani~\cite{mir4}.

It is natural to consider the asymptotic behaviour of $V_{g,n}(N\x)$ for a fixed $\x = (x_1, x_2, \ldots, x_n)$ as $N$ approaches infinity. We are thus motivated to analyse a certain normalisation of the Weil--Petersson form on ${\mathcal M}_{g,n}(N\x)$ as $N$ approaches infinity~\cite{do2}. In the limit, we obtain a 2-form originally defined by Kontsevich in his proof of Witten's conjecture~\cite{kon}. In this way, we obtain yet another proof of the Witten--Kontsevich theorem which makes explicit the connection between the work of Kontsevich and Mirzakhani.

The structure of the article is as follows.
\begin{itemize} [topsep=0pt,itemsep=1pt]
\item In Section~\ref{moduli-spaces}, we introduce hyperbolic surfaces and their moduli spaces via Teichm\"{u}ller theory. We discuss the Weil--Petersson symplectic structure of the space ${\mathcal M}_{g,n}(\LL)$ and its Deligne--Mumford compactification. We conclude by showing how the hyperbolic geometry leads to a cell decomposition of the moduli space based on the combinatorial notion of a ribbon graph.

\item In Section~\ref{WP-volumes}, we begin the study of Weil--Petersson volumes. Some early results in the area are presented, followed by some preparatory remarks on symplectic reduction. We then show how Mirzakhani applies this technique to prove the polynomiality of Weil--Petersson volumes~\cite{mir1}.

\item In Section~\ref{recursive-formula}, we discuss Mirzakhani's recursion for Weil--Petersson volumes~\cite{mir2}. We calculate the volume of ${\mathcal M}_{1,1}(0)$ as a motivating example. The generalised McShane identity is presented and then used as one of the main ingredients in the proof of Mirzakhani's recursive formula. We consider some applications of the recursion, including a proof of the Witten--Kontsevich theorem.

\item In Section~\ref{WP-limits}, we analyse Weil--Petersson volumes under various limits. In particular, we examine the recent results on the behaviour of $V_{g,n}(\LL)$ as one of the lengths approaches $2\pi i$~\cite{do-nor}, in the large $g$ limit~\cite{mir4}, and as the lengths approach infinity~\cite{do2}.
\end{itemize}

\section{Moduli spaces of hyperbolic surfaces} \label{moduli-spaces}

\subsection{Hyperbolic surfaces}

Consider the smooth surface $\Sigma_{g,n}$ with genus $g$ and $n$ boundary components. A {\em hyperbolic surface of type $(g,n)$} is the surface $\Sigma_{g,n}$ equipped with a complete Riemannian metric of constant curvature $-1$. We restrict our attention to the case when every boundary component is smooth and totally geodesic. A mild restriction on the pair of non-negative integers $(g,n)$ is required due to the Gauss-Bonnet theorem. In its simplest form, it states that the integral of the Gaussian curvature over a surface with Riemannian metric and totally geodesic boundary is equal to $2\pi$ multiplied by the Euler characteristic of the surface.
\[
\int_S K\,dA = 2\pi \chi(S)
\]
So for a metric of constant curvature $-1$ to exist, the Euler characteristic must be negative --- in other words, $2 - 2g - n < 0$. This is a rather mild restriction since it only prohibits the pairs $(0, 0)$, $(0, 1)$, $(0,2)$ and $(1,0)$. Note that these exceptional cases are precisely the pairs $(g, n)$ for which a Riemann surface of genus $g$ with $n$ punctures possesses infinitely many automorphisms.

An alternative definition of a hyperbolic surface uses the notion of an atlas. Thus, we define a hyperbolic surface to be a smooth surface covered by charts $\phi_U: U \to \mathbb{H}^2$ which map open subsets of the surface homeomorphically onto their image in the hyperbolic plane. We require that if $U \cap V \neq \emptyset$, then the two charts are compatible in the sense that the transition function $\phi_V \circ \phi_U^{-1}: \phi_U(U \cap V) \to \phi_V(U \cap V)$ is an isometry. As usual, a hyperbolic surface is defined by a maximal atlas --- that is, a maximal collection of compatible charts. The two definitions provided coincide since any two-dimensional domain with a Riemannian metric of constant curvature $-1$ is locally isometric to a subset of the hyperbolic plane. Which to use as a working definition is largely a matter of taste, though it is advantageous to keep both viewpoints in mind.

Yet another way to define a hyperbolic surface is via its universal cover. Every hyperbolic surface $S$ with geodesic boundary has a universal cover isometric to a convex domain in $\mathbb{H}^2$ with geodesic boundary. Therefore, hyperbolic surfaces with geodesic boundary arise from taking the quotient of a convex domain in $\mathbb{H}^2$ with geodesic boundary by a discrete subgroup of $PSL(2, \mathbb{R})$, the group of orientation-preserving isometries in the hyperbolic plane. The following result gives one of the important properties of hyperbolic surfaces.

\begin{proposition} \label{geodesic}
On a hyperbolic surface, non-trivial homotopy classes of closed curves have unique geodesic representatives. Furthermore, such geodesic representatives realise minimal intersection and self-intersection numbers. In particular, every simple closed curve is homotopic to a simple closed geodesic.
\end{proposition}

For an $n$-tuple $\LL = (L_1, L_2, \ldots, L_n)$ of positive real numbers, we define the {\em moduli space of hyperbolic surfaces} as follows.
\[
{\mathcal M}_{g,n}(\LL) = \left. \left\{ (S, \beta_1, \beta_2, \ldots, \beta_n) \; \middle \vert \; \begin{array}{l} S \text{ is a hyperbolic surface of type } (g,n) \text{ with boundary} \\ \text{components } \beta_1, \beta_2, \ldots, \beta_n \text{ of lengths } L_1, L_2, \ldots, L_n \end{array} \right\} \right/ \sim
\]
Here, $(S, \beta_1, \beta_2, \ldots, \beta_n) \sim (T, \gamma_1, \gamma_2, \ldots, \gamma_n)$ if and only if there exists an isometry from $S$ to $T$ which sends $\beta_k$ to $\gamma_k$ for all $k$.

Note that when the length of a boundary component approaches zero, we obtain a hyperbolic cusp in the limit. By a {\em hyperbolic cusp}, we mean a subset of the surface isometric to $\mathbb{R} / \mathbb{Z} \times [2, \infty)$, in the Poincar\'{e} upper half-plane model for $\mathbb{H}^2$. In particular, we denote the moduli space of genus $g$ hyperbolic surfaces with $n$ labelled cusps by ${\mathcal M}_{g,n}(\mathbf{0})$.

One of the goals of this article is to explain how the study of moduli spaces of hyperbolic surfaces can lead to results concerning moduli spaces of curves. This is possible through the interplay between hyperbolic surfaces, Riemann surfaces, and algebraic curves. It is well-known that the category of smooth algebraic curves is equivalent to the category of compact Riemann surfaces. Due to this equivalence, the boundary between these two fields is rather porous, with techniques from complex analysis flowing into algebraic geometry and vice versa. The uniformisation theorem then provides a one-to-one correspondence between Riemann surfaces and hyperbolic surfaces.

\begin{theorem}[Uniformisation theorem] \label{uniformisation}
Every Riemannian metric on a surface is conformally equivalent to a complete constant curvature metric. If the Euler characteristic of the surface is negative and we require that the curvature is $-1$, then the metric is unique.
\end{theorem}

From the previous discussion, a smooth genus $g$ algebraic curve with $n$ labelled points corresponds to a genus $g$ Riemann surface with $n$ labelled points, which we usually think of as punctures.
\[
\left\{ \begin{array}{c}
\text{smooth algebraic curves with} \\
\text{genus $g$ and $n$ marked points} \end{array} \right\} 
\longleftrightarrow 
\left\{ \begin{array}{c}
\text{Riemann surfaces with} \\
\text{genus $g$ and $n$ punctures} \end{array} \right\}
\]
The complex structure on a Riemann surface defines a conformal class of metrics which, by the uniformisation theorem, contains a hyperbolic metric as long as $2-2g-n < 0$. Furthermore, if we demand that the resulting surface is complete, then this hyperbolic metric is unique and gives each puncture the structure of a hyperbolic cusp. So we have the following one-to-one correspondence.
\[
\left\{ \begin{array}{c}
\text{Riemann surfaces with} \\
\text{genus $g$ and $n$ punctures} \end{array} \right\}
\longleftrightarrow 
\left\{ \begin{array}{c}
\text{hyperbolic surfaces with} \\
\text{genus $g$ and $n$ cusps} \end{array} \right\}
\]

At the moment, the moduli space of hyperbolic surfaces ${\mathcal M}_{g,n}(\LL)$ has only been defined as a set. In Section~\ref{teichmuller-theory}, we will see that it possesses not only a topology, but also an orbifold structure. The equivalence above sets up a bijection between the moduli space of curves ${\mathcal M}_{g,n}$ and the moduli space of hyperbolic surfaces ${\mathcal M}_{g,n}(\mathbf{0})$. This map respects the topology of both spaces as well as the structure-preserving automorphism groups. As a result, ${\mathcal M}_{g,n}$ and ${\mathcal M}_{g,n}(\mathbf{0})$ are homeomorphic as orbifolds.

\subsection{Teichm\"{u}ller theory} \label{teichmuller-theory}

Teichm\"{u}ller theory will enable us to construct the moduli space of hyperbolic surfaces ${\mathcal M}_{g,n}(\LL)$ and endow it with a natural symplectic structure. Begin by fixing a smooth surface $\Sigma_{g,n}$ with genus $g$ and $n$ boundary components labelled from 1 up to $n$, where $2-2g-n < 0$. Define a {\em marked hyperbolic surface of type $(g,n)$} to be a pair $(S, f)$ where $S$ is a hyperbolic surface and $f: \Sigma_{g,n} \to S$ is a diffeomorphism. We call $f$ the marking of the hyperbolic surface and define the {\em Teichm\"{u}ller space} as follows.
\[
{\mathcal T}_{g,n}(\LL) = \left. \left\{ (S, f) \; \middle \vert \; \begin{array}{l} (S, f) \text{ is a marked hyperbolic surface of type } (g,n) \text{ with} \\ \text{boundary components of lengths } L_1, L_2, \ldots, L_n \end{array} \right\} \right/ \sim
\]
Here, $(S, f) \sim (T, g)$ if and only if there exists an isometry $\phi: S \to T$ such that $\phi \circ f$ is isotopic to $g$.

One can informally think of Teichm\"{u}ller space as the space of deformations of the hyperbolic structure on a given hyperbolic surface. For example, consider applying a hyperbolic Dehn twist to a marked hyperbolic surface. By this, we mean cutting along a simple closed geodesic, twisting the two sides relative to each other, and gluing the two sides back together. This gives a one parameter family of deformations of the hyperbolic structure. Once a full twist has been applied, the end result is a hyperbolic surface isometric to the original and hence, corresponds to the same point in the moduli space. On the other hand, the end result has a different marking to the original and hence, corresponds to a different point in the Teichm\"{u}ller space.

The geometry and topology of Teichm\"{u}ller space will become much more apparent once we define global coordinates, known as Fenchel--Nielsen coordinates. We use the idea that pairs of pants --- spheres with three boundary components --- can be used as building blocks to create surfaces with negative Euler characteristic. Start by considering a pants decomposition of the surface $\Sigma_{g,n}$, which is a collection of disjoint simple closed curves whose complement is a disjoint union of pairs of pants. Alternatively, a pants decomposition is a maximal collection of disjoint simple closed curves such that no curve is parallel to the boundary and no two are homotopic. Since the Euler characteristic is additive over surfaces glued along circles, the number of pairs of pants in any such decomposition must be $-\chi(\Sigma_{g,n}) = 2g-2+n$. Some simple combinatorics can be used to show that every pants decomposition of $\Sigma_{g,n}$ consists of precisely $3g-3+n$ simple closed curves.

A marking $f: \Sigma_{g,n} \to S$ maps a fixed pants decomposition to a collection of simple closed curves on $S$, each of which is homotopic to a unique simple closed geodesic by Proposition~\ref{geodesic}. Denote these simple closed geodesics by $\gamma_1, \gamma_2, \ldots, \gamma_{3g-3+n}$ and let their lengths be $\ell_1, \ell_2, \ldots, \ell_{3g-3+n}$, respectively. Cutting $S$ along $\gamma_1, \gamma_2, \ldots, \gamma_{3g-3+n}$ leaves a disjoint union of $2g-2+n$ hyperbolic pairs of pants. The following elementary result guarantees that the lengths $\ell_1, \ell_2, \ldots, \ell_{3g-3+n}$ are sufficient to determine the hyperbolic structure on each pair of pants.

\begin{lemma} \label{pair-of-pants}
There exists a unique hyperbolic pair of pants up to isometry with geodesic boundary components of prescribed non-negative length. As usual, we interpret a geodesic boundary component with length zero as a hyperbolic cusp. The three simple geodesic arcs perpendicular to the boundary components and joining them in pairs are referred to as seams. Cutting along the seams decomposes a hyperbolic pair of pants into two congruent right-angled hexagons.
\end{lemma}

Note that the lengths $\ell_1, \ell_2, \ldots, \ell_{3g-3+n}$ provide insufficient information to reconstruct the hyperbolic structure on all of $S$, since there are infinitely many ways to glue together the pairs of pants. This extra gluing information is stored in the twist parameters, which we denote by $\tau_1, \tau_2, \ldots, \tau_{3g-3+n}$. To construct them, fix a collection $C$ of disjoint smooth curves on $\Sigma_{g,n}$ which are either closed or have endpoints on the boundary. We require that $C$ meets the pants decomposition transversely, such that its restriction to any pair of pants consists of three disjoint arcs, connecting the boundary components pairwise. Now to construct the twist parameter $\tau_k$, take a curve $\gamma \in C$ such that $f(\gamma)$ meets $\gamma_k$. Homotopic to $f(\gamma)$, relative to the boundary of $S$, is a unique length-minimising piecewise geodesic curve which is entirely contained in the seams of the pairs of pants and the curves $\gamma_1, \gamma_2, \ldots, \gamma_{3g-3+n}$. The twist parameter $\tau_k$ is the signed distance that this curve travels along $\gamma_k$, according to the following sign convention. Lemma~\ref{pair-of-pants} guarantees that the twist parameter is independent of the choice of curve $\gamma \in C$.

\vspace{6pt}
\begin{center}
\begin{overpic}[unit=1mm]{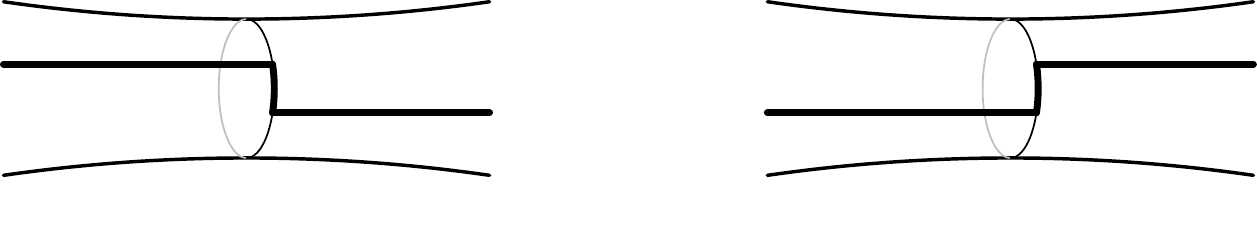}
\put(4.5,3){negative twist parameter}
\put(66,3){positive twist parameter}
\end{overpic}
\end{center}
\vspace{-6pt}

For further details, one can consult Thurston's book~\cite{thu}, in which he writes the following.
\begin{quote}
``That a twist parameter takes values in $\mathbb{R}$, rather than $S^1$, tends to be a confusing issue \ldots But, remember, to determine a point in Teichm\"{u}ller space we need to consider how many times the leg of the pajama suit is twisted before it fits onto the baby's foot.''
\end{quote}

More prosaically, the length parameters and the twist parameters modulo the length parameters are sufficient to reconstruct the hyperbolic structure on $S$. However, to recover the marking as well, it is necessary to consider the twist parameters as elements of $\mathbb{R}$. Despite the fact that the length and twist parameters --- collectively known as {\em Fenchel--Nielsen coordinates} --- depend on the choice of pants decomposition and the construction of twist parameters, we have the following result.

\begin{theorem}
The map ${\mathcal T}_{g,n}(\LL) \to \mathbb{R}_+^{3g-3+n} \times \mathbb{R}^{3g-3+n}$, which associates to a marked hyperbolic surface its length and twist parameters, is a bijection. In fact, if Teichm\"{u}ller space is considered with its natural topology, then the map is a homeomorphism.
\end{theorem}

Clearly, there is a projection map ${\mathcal T}_{g,n}(\LL) \to {\mathcal M}_{g,n}(\LL)$ which simply forgets the marking. In fact, the moduli space is obtained as a quotient of Teichm\"{u}ller space by a group action. Define the {\em mapping class group} as
\[
\textup{Mod}_{g,n} = \textup{Diff}^+(\Sigma_{g,n})/\textup{Diff}^+_0(\Sigma_{g,n}),
\]
where $\textup{Diff}^+$ denotes the group of orientation preserving diffeomorphisms fixing the boundary components and $\textup{Diff}^+_0$ denotes the normal subgroup consisting of those diffeomorphisms isotopic to the identity. There is a natural action of the mapping class group on Teichm\"{u}ller space such that $[\phi] \in \textup{Mod}_{g,n}$ sends the marked hyperbolic surface $(X, f)$ to the marked hyperbolic surface $(X, f \circ \phi)$. The moduli space ${\mathcal M}_{g,n}(\LL)$ is obtained by taking the quotient of the Teichm\"{u}ller space ${\mathcal T}_{g,n}(\LL)$ by the action of the mapping class group $\textup{Mod}_{g,n}$.

\begin{proposition}
The action of $\textup{Mod}_{g,n}$ on ${\mathcal T}_{g,n}(\LL)$ is properly discontinuous, though not necessarily free. Therefore, the quotient ${\mathcal M}_{g,n}(\LL) = \T_{g,n}(\LL) / \textup{Mod}_{g,n}$ is an orbifold of dimension $6g-6+2n$.
\end{proposition}

So the moduli space of hyperbolic surfaces ${\mathcal M}_{g,n}(\LL)$ has not only a topology, but also an orbifold structure. The orbifold group at a point is canonically isomorphic to the automorphism group of the corresponding hyperbolic surface. However, the situation is not so bad, since the following theorem --- which follows from results of Boggi and Pikaart~\cite{bog-pik} --- allows one to make sense of calculations on the orbifold by lifting to a finite cover. For this reason, it is convenient for us to consider the cohomology of moduli spaces with rational, rather than integral, coefficients.

\begin{theorem} \label{manifold-cover}
A finite cover $\widetilde{\mathcal M}_{g,n}(\LL) \to {\mathcal M}_{g,n}(\LL)$ exists such that $\widetilde{\mathcal M}_{g,n}(\LL)$ is a smooth manifold.
\end{theorem}

We have only very briefly touched upon the vast area that is Teichm\"{u}ller theory. For more information, see the {\em Handbook of Teichm\"{u}ller theory}~\cite{pap1, pap2}.

\subsection{Symplectification and compactification} \label{symplectification-compactification}

The Teichm\"{u}ller space $\T_{g,n}(\LL)$ can be endowed with the canonical symplectic form
\[
\omega = \sum_{k=1}^{3g-3+n} d\ell_k \wedge d\tau_k
\]
using the Fenchel--Nielsen coordinates. Although this is a rather trivial statement, it is a deep fact that this form is invariant under the action of the mapping class group. Therefore, $\omega$ descends to a symplectic form on the quotient, namely the moduli space ${\mathcal M}_{g,n}(\LL)$. This is referred to as the {\em Weil--Petersson symplectic form} and we will also denote it by $\omega$. Its existence allows us to introduce the techniques of symplectic geometry to the study of moduli spaces. For all values of $\LL$, the spaces ${\mathcal M}_{g,n}(\LL)$ are diffeomorphic to each other, but not necessarily symplectomorphic to each other. It is therefore natural to ask how the symplectic structure varies as $\LL$ varies, a question which we will pursue in Section~\ref{polynomiality}.

The uniformisation theorem allows us to deduce that the moduli spaces ${\mathcal M}_{g,n}(\LL)$ are not only diffeomorphic to each other, but also diffeomorphic to the moduli space of curves ${\mathcal M}_{g,n}$. It is often more natural to work with the Deligne--Mumford compactification $\M_{g,n}$, obtained by introducing the notion of a stable algebraic curve --- see Appendix~\ref{intersection-theory} for the relevant definitions. There is an analogous construction in the hyperbolic setting, where a node of an algebraic curve corresponds to degenerating the length of a simple closed curve on a hyperbolic surface to zero. This intuition leads to the following construction of $\overline{\mathcal T}_{g,n}(\LL)$, the Teichm\"{u}ller space of marked stable hyperbolic surfaces. Define a {\em stable hyperbolic surface of type $(g,n)$} to be a pair $(S, M)$ where $S$ is a surface of genus $g$ with $n$ labelled boundary components and $M$ is a collection of disjoint simple closed curves on $S$, none of which are homotopic to the boundary. We require that $S \setminus M$ be endowed with a finite area hyperbolic metric such that the boundary components are geodesic. It is useful to think of a stable hyperbolic surface as a collection of hyperbolic surfaces whose cusps have been formally identified in pairs. As usual, we refer to a diffeomorphism $f: \Sigma_{g,n} \to S$ as a marking and define the {\em compactified Teichm\"{u}ller space} as follows.
\[
\overline{\mathcal T}_{g,n}(\LL) = \left. \left\{ (S, M, f) \; \middle \vert \; \begin{array}{l} (S, M, f) \text{ is a marked stable hyperbolic surface of type } (g,n) \\ \text{with boundary components of lengths } L_1, L_2, \ldots, L_n \end{array} \right\} \right/ \sim
\]
Here, $(S, M, f) \sim (T, N, g)$ if and only if there exists a homeomorphism $\phi: S \to T$ such that $\phi(M) = N$, $\phi$ restricted to $S \setminus M$ is an isometry, and $\phi \circ f$ is isotopic to $g$ on each connected component of $\Sigma_{g,n} \setminus f^{-1}(M)$. Once again, the mapping class group acts on the compactified Teichm\"{u}ller space and one may define the {\em Deligne--Mumford compactification} of the moduli space of hyperbolic surfaces as
\[
\M_{g,n}(\LL) = \overline{\mathcal T}_{g,n}(\LL) / \textup{Mod}_{g,n}.
\]

Let us make some remarks on the Deligne--Mumford compactification. First, the compactification locus $\M_{g,n}(\LL) \setminus {\mathcal M}_{g,n}(\LL)$ is a union of submanifolds of positive codimension. Second, $\M_{g,n}(\mathbf{0})$ can be canonically identified with the Deligne--Mumford compactification of the moduli space of curves $\M_{g,n}$ via the uniformisation theorem. Hence, $\M_{g,n}(\mathbf{0})$ possesses a natural complex structure, whereas the moduli space $\overline{\mathcal M}_{g,n}(\LL)$ does not, for $\LL \neq \mathbf{0}$. However, by the work of Wolpert, the Fenchel--Nielsen coordinates do induce a real analytic structure~\cite{abi, wol1}.

Wolpert~\cite{wol1} used the real analytic structure on ${\mathcal M}_{g,n}(\LL)$ to show that the Weil--Petersson form extends smoothly to a closed non-degenerate form on the Deligne--Mumford compactification $\M_{g,n}(\LL)$. In the particular case $\LL = \mathbf{0}$, he showed that this extension defines a cohomology class $[\omega] \in H^2(\M_{g,n}, \mathbb{R})$ which satisfies the following.

\begin{theorem} \label{kappa1}
The de Rham cohomology class of the Weil--Petersson symplectic form on $\overline{\mathcal M}_{g,n}(\mathbf{0})$ satisfies $[\omega] = 2\pi^2 \kappa_1 \in H^2(\M_{g,n}, \mathbb{R})$, where $\kappa_1$ denotes the first Mumford--Morita--Miller class.
\end{theorem}

\subsection{Combinatorial moduli space} \label{combinatorial-moduli-space}

An important notion in the study of moduli spaces of curves is the combinatorial structure known in the literature as a ribbon graph or fatgraph. A {\em ribbon graph of type $(g,n)$} is essentially the 1-skeleton of a cell decomposition of a genus $g$ surface with $n$ faces. We require the vertices to have degree at least three and the faces to be labelled from 1 up to $n$. Note that such a graph may possibly have loops or multiple edges. The orientation of the surface gives a cyclic ordering to the oriented edges pointing toward each vertex. Conversely, given the underlying graph and the cyclic ordering of the oriented edges pointing toward each vertex, the genus of the surface and its cell decomposition may be recovered. This is accomplished by using the extra structure to thicken the graph into a surface with boundaries. These boundaries may then be filled in with disks to produce a closed surface with an associated cell decomposition.

One usually draws ribbon graphs with the convention that the cyclic ordering of the oriented edges pointing toward each vertex is induced by the orientation of the page. For example, the following diagram shows a ribbon graph of type $(1,1)$ as well as the surface obtained by thickening the graph.

\begin{center}
\includegraphics{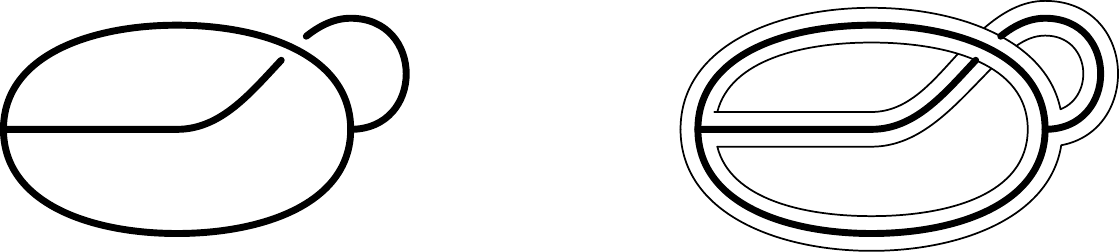}
\end{center}

It is often useful to think of a ribbon graph in the following more precise way. Given a cell decomposition $\Gamma$ of a surface, let $X$ denote the set of its oriented edges and let $s_0$ be the permutation on $X$ which cyclically permutes all oriented edges pointing toward the same vertex in an anticlockwise manner. Also, let $s_1$ be the permutation on $X$ which interchanges each pair of oriented edges which correspond to the same underlying edge. The set $X / \langle s_0 \rangle$ is canonically equivalent to the set of vertices of $\Gamma$ while the set $X / \langle s_1 \rangle$ is canonically equivalent to the set of edges of $\Gamma$. Furthermore, if we let $s_2 = s_1 s_0^{-1}$, then the set $X / \langle s_2 \rangle$ is canonically equivalent to the set of faces of $\Gamma$. Therefore, one can alternatively define a ribbon graph to be a triple $(X, s_0, s_1)$ where $X$ is a finite set, $s_0$ is a permutation on $X$ without fixed points or transpositions, and $s_1$ is an involution on $X$ without fixed points. We also require a labelling in the form of a bijection from $X / \langle s_2 \rangle$ to $\{1, 2, \ldots, n\}$. Define two ribbon graphs $(X, s_0, s_1)$ and $(\overline{X}, \overline{s}_0, \overline{s}_1)$ to be isomorphic if and only if there exists a bijection $f: X \to \overline{X}$ such that $f \circ s_0 = \overline{s}_0 \circ f$ and $f \circ s_1 = \overline{s}_1 \circ f$. We also impose the condition that $f$ must preserve the labelling of the boundary components. A ribbon graph automorphism is, of course, an isomorphism from a ribbon graph to itself. The set of automorphisms of a ribbon graph $\Gamma$ forms a group which is denoted by $\mathrm{Aut}(\Gamma)$.

A ribbon graph with a positive real number assigned to each edge is referred to as a {\em metric ribbon graph}. The metric associates to each face in the cell decomposition a perimeter, which is simply the sum of the numbers appearing around the boundary of the face. We define the {\em combinatorial moduli space} as follows.
\[
{\mathcal MRG}_{g,n}(\LL) = \left. \left\{ \begin{array}{c} \text{metric ribbon graphs of type } (g,n) \\ \text{with perimeters } L_1, L_2, \ldots, L_n \end{array} \right\} \right/ \sim
\]
Here, two metric ribbon graphs are equivalent if and only if there exists an isometry between them which corresponds to a ribbon graph automorphism.

For a ribbon graph $\Gamma$ of type $(g, n)$, consider ${\mathcal MRG}_\Gamma(\LL) \subseteq {\mathcal MRG}_{g,n}(\LL)$, the subset consisting of those metric ribbon graphs whose underlying ribbon graph is $\Gamma$. Note that ${\mathcal MRG}_\Gamma(\LL)$ can be naturally identified with the following quotient of a possibly empty polytope by a finite group.
\[
{\mathcal MRG}_\Gamma(\LL) \cong \left. \left\{ \mathbf{e} \in \mathbb{R}_+^{E(\Gamma)} \; \middle \vert \; A_\Gamma \mathbf{e} = \LL \right\} \right/ \mathrm{Aut}(\Gamma)
\]
Here, $\mathbf{e}$ represents the lengths of the edges in the metric ribbon graph, $E(\Gamma)$ denotes the edge set of $\Gamma$, and $A_\Gamma$ is the linear map which represents the adjacency between faces and edges in the cell decomposition corresponding to $\Gamma$. Thus, ${\mathcal MRG}_\Gamma(\LL)$ is an orbifold cell and these naturally glue together via edge degenerations --- in other words, when an edge length goes to zero, the edge contracts to give a ribbon graph with fewer edges. So this cell decomposition for ${\mathcal MRG}_{g,n}(\LL)$ equips it with not only a topology, but also an orbifold structure. The main reason for considering ${\mathcal MRG}_{g,n}(\LL)$ is the following result.

\begin{theorem} \label{bowditch-epstein}
The moduli spaces ${\mathcal M}_{g,n}(\LL)$ and ${\mathcal MRG}_{g,n}(\LL)$ are homeomorphic as orbifolds.
\end{theorem}

One can prove this fact by generalising the work of Bowditch and Epstein~\cite{bow-eps}, who consider the case of cusped hyperbolic surfaces. The main idea is to associate to a hyperbolic surface $S$ with geodesic boundary its spine $\Gamma(S)$. For every point $p \in S$, let $n(p)$ denote the number of shortest paths from $p$ to the boundary. Generically, we have $n(p) = 1$ and we define the {\em spine} as
\[
\Gamma(S) = \{p \in S \mid n(p) \geq 2\}.
\]
The locus of points with $n(p) = 2$ consists of a disjoint union of open geodesic segments. These correspond precisely to the edges of a graph embedded in $S$. The locus of points with $n(p) \geq 3$ forms a finite set which corresponds to the set of vertices of the aforementioned graph. In fact, if $n(p) \geq 3$, then the corresponding vertex will have degree $n(p)$. In this way, $\Gamma(S)$ has the structure of a ribbon graph. Furthermore, it is a deformation retract of the original hyperbolic surface, so if $S$ is a hyperbolic surface of type $(g,n)$, then $\Gamma(S)$ will be a ribbon graph of type $(g,n)$.

Now for each vertex $p$ of $\Gamma(S)$, consider the $n(p)$ shortest paths from $p$ to the boundary. We refer to these geodesic segments as {\em ribs} and note that they are perpendicular to the boundary of $S$. The diagram below shows part of a hyperbolic surface, along with its spine and ribs. Cutting $S$ along its ribs leaves a collection of hexagons, each with four right angles and a reflective axis of symmetry along one of the diagonals. In fact, this diagonal is one of the edges of $\Gamma(S)$ and we assign to it the length of the side of the hexagon which lies along the boundary of $S$. Of course, there are two such sides --- however, the reflective symmetry guarantees that they are equal in length. In this way, $\Gamma(S)$ becomes a metric ribbon graph of type $(g,n)$. By construction, the perimeters of $\Gamma(S)$ correspond precisely with the lengths of the boundary components of $S$, so we have a map $\Gamma: {\mathcal M}_{g,n}(\LL) \to {\mathcal MRG}_{g,n}(\LL)$. It is possible, though more difficult, to construct the inverse map $S: {\mathcal MRG}_{g,n}(\LL) \to {\mathcal M}_{g,n}(\LL)$ and show that it preserves the orbifold structure of the moduli spaces. The omitted details may be found elsewhere in the literature~\cite{bow-eps, do1}.

\begin{center}
\includegraphics{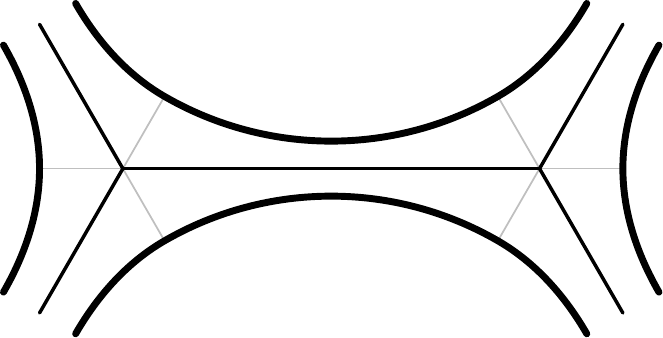}
\end{center}
\vspace{-6pt}

The notion of the combinatorial moduli space is crucial to Kontsevich's proof of Witten's conjecture concerning intersection numbers on $\M_{g,n}$~\cite{kon}. We remark that Kontsevich uses a different construction of the combinatorial moduli space which begins with punctured Riemann surfaces. Metric ribbon graphs arise via the existence of Jenkins--Strebel quadratic differentials on Riemann surfaces, an observation which Kontsevich attributes to Harer, Mumford, Penner and Thurston.

\section{Weil--Petersson volumes} \label{WP-volumes}

\subsection{Early results} \label{early-results}

By raising the Weil--Petersson symplectic form to the appropriate exterior power, one obtains the following volume form on the Teichm\"{u}ller space ${\mathcal T}_{g,n}(\LL)$.
\[
\frac{\omega^{3g-3+n}}{(3g-3+n)!} = d\ell_1 \wedge d\tau_1 \wedge d\ell_2 \wedge d\tau_2 \wedge \cdots \wedge d\ell_{3g-3+n} \wedge d\tau_{3g-3+n}
\]
Of course, ${\mathcal T}_{g,n}(\LL)$ has infinite volume with respect to this form. However, the action of the mapping class group is such that the volume of ${\mathcal M}_{g,n}(\LL)$ is finite. One way to see this is via Wolpert's observation~\cite{wol1} that the Weil--Petersson symplectic form extends smoothly to a closed non-degenerate form on the Deligne--Mumford compactification $\overline{\mathcal M}_{g,n}(\LL)$. Thus, let us define the {\em Weil--Petersson volume}
\[
V_{g,n}(\LL) = \int_{{\mathcal M}_{g,n}(\LL)} \frac{\omega^{3g-3+n}}{(3g-3+n)!}.
\]
Note that we may instead choose to integrate over the Deligne--Mumford compactification, since the compactification locus is the union of submanifolds of positive codimension.

Below we present a brief selection of some early results concerning Weil--Petersson volumes. When comparing these results with those in the literature, there may be some discrepancy due to two issues. First, there are distinct normalisations of the Weil--Petersson symplectic form which differ by a factor of two. We have scaled the results, where appropriate, to correspond to the Weil--Petersson symplectic form defined in Section~\ref{symplectification-compactification}. Second, one must treat the special cases of $V_{1,1}(L_1)$ and $V_{2,0}$ with some care. This is due to the fact that every point on ${\mathcal M}_{1,1}(L)$ and ${\mathcal M}_{2,0}$ is an orbifold point, generically with orbifold group equal to $\mathbb{Z} / 2\mathbb{Z}$. Cleaner statements of results are obtained if one considers $V_{1,1}(L_1)$ and $V_{2,0}$ as orbifold volumes --- in this case, half of the true volumes. The upshot is that one should not be alarmed if results concerning Weil--Petersson volumes from distinct sources differ by a factor which is a power of two.

\begin{itemize} [topsep=0pt,itemsep=1pt]
\item Wolpert~\cite{wol2, wol3} proved that $V_{0,4}(0, 0, 0, 0) = 2\pi^2$, $V_{1,1}(0) = \frac{\pi^2}{12}$ and $V_{g,n}(\mathbf{0}) = q (2\pi^2)^{3g-3+n}$ for some rational number $q$. This last fact is a corollary of Theorem~\ref{kappa1}, from which it follows that $q = \int_{\M_{g,n}} \kappa_1^{3g-3+n}$.
\item Penner~\cite{pen} proved that $V_{1,2}(0, 0) = \frac{\pi^4}{4}$.
\item Zograf~\cite{zog1} proved that $V_{0,n}(\mathbf{0}) = \frac{(2\pi^2)^{n-3}}{(n-3)!} a_n$, where $a_3 = 1$ and
\[
a_n = \frac{1}{2} \sum_{k=1}^{n-3} \frac{k(n-k-2)}{n-1} \binom{n-4}{k-1} \binom{n}{k+1} a_{k+2} a_{n-k} \qquad \text{for } n \geq 4.
\]
\item N\"{a}\"{a}t\"{a}nen and Nakanishi~\cite{nak-naa} calculated the Weil--Petersson volumes of the two-dimensional moduli spaces.
\begin{align*}
V_{0,4}(L_1, L_2, L_3, L_4) &= \frac{1}{2}(L_1^2 + L_2^2 + L_3^2 + L_4^2 + 4\pi^2) \\
V_{1,1}(L_1) &= \frac{1}{48}(L_1^2 + 4\pi^2)
\end{align*}
\end{itemize}

N\"{a}\"{a}t\"{a}nen and Nakanishi's result shows that $V_{1,1}(L_1)$ and $V_{0,4}(L_1, L_2, L_3, L_4)$ are both polynomials in the squares of the boundary lengths. That this is the case for all Weil--Petersson volumes $V_{g,n}(\LL)$ was proven by Mirzakhani in two distinct ways. In Section~\ref{polynomiality}, we discuss the first of Mirzakhani's proofs, which uses symplectic reduction in a fundamental way~\cite{mir1}.

\subsection{Symplectic reduction}

Symplectic geometry has its origins in the mathematical formulation and generalisation of the phase space of a classical mechanical system. Physicists have often taken advantage of the fact that when a symmetry group of dimension $n$ acts on a system, then the number of degrees of freedom for the positions and momenta can be reduced by $2n$. The analogous mathematical phenomenon is known as {\em symplectic reduction}. More precisely, take a symplectic manifold $(M, \omega)$ of dimension $2d$ with a
\[
T^n = \underbrace{S^1 \times S^1 \times \ldots \times S^1}_{n \text{ times}}
\]
action that preserves the symplectic form. Furthermore, suppose that this action is the Hamiltonian flow for the moment map $\mu: M \to \mathbb{R}^n$ and that $\mathbf{0}$ is a regular value of $\mu$. By this we mean that the Hamiltonian vector field $X_\mu$ defined by the equation $\omega(X_\mu, \cdot) = dH(\cdot)$ generates the action. Since $T^n$ must act on the level sets of $\mu$, we can define $M_{\mathbf{0}} = \mu^{-1}(\mathbf{0}) / T^n$.

\begin{theorem}[Marsden--Weinstein theorem]
The orbit space $M_{\mathbf{0}} = \mu^{-1}(\mathbf{0}) / T^n$ is a symplectic manifold of dimension $2d - 2n$ with respect to the unique 2-form $\omega_{\mathbf{0}}$ which satisfies $i^*\omega = \pi^*\omega_{\mathbf{0}}$. Here, $\pi: \mu^{-1}(\mathbf{0}) \to M_{\mathbf{0}}$ and $i: \mu^{-1}(\mathbf{0}) \to M$ are the natural projection and inclusion maps. 
\end{theorem}

Since $\mathbf{0}$ is a regular value, there exists an $\varepsilon > 0$ such that all $\mathbf{a} \in \mathbb{R}^n$ satisfying $|\mathbf{a}| < \varepsilon$ are also regular values. So it is possible to define symplectic manifolds $(M_{\mathbf{a}}, \omega_{\mathbf{a}})$ for all such $\mathbf{a}$. If we think of the $T^n$ action as $n$ commuting circle actions, then the $k$th copy of $S^1$ induces a circle bundle ${\mathcal S}_k$ on $M_{\mathbf{0}}$. The variation of the symplectic form $\omega_{\mathbf{0}}$ can be described in terms of the first Chern classes $\phi_k = c_1({\mathcal S}_k)$.

\begin{theorem} \label{formvar}
For $\mathbf{a} = (a_1, a_2, \ldots, a_n)$ sufficiently close to $\mathbf{0}$, $(M_{\mathbf{a}}, \omega_{\mathbf{a}})$ is symplectomorphic to $M_{\mathbf{0}}$ equipped with a symplectic form whose cohomology class is equal to $[\omega_{\mathbf{0}}] + a_1 \phi_1 + a_2 \phi_2 + \cdots + a_n \phi_n$.
\end{theorem}

From this theorem, one obtains as a direct corollary an expression for the variation of the volume.

\begin{corollary} \label{volumevar}
For $\mathbf{a} = (a_1, a_2, \ldots, a_n)$ sufficiently close to $\mathbf{0}$, the volume of $(M_{\mathbf{a}}, \omega_{\mathbf{a}})$ is a polynomial in $a_1, a_2, \ldots, a_n$ of degree $d = \frac{1}{2} \dim(M_{\mathbf{a}})$ given by the formula
\[
\sum_{|\bm{\alpha}| + m = d} \frac{\int_{M_{\mathbf{0}}} \phi_1^{\alpha_1} \phi_2^{\alpha_2} \cdots \phi_n^{\alpha_n} \omega^{m}}{\alpha_1! \alpha_2! \cdots \alpha_n! m!} a_1^{\alpha_1} a_2^{\alpha_2} \cdots a_n^{\alpha_n}.
\]
\end{corollary}

For an introduction to symplectic geometry, we recommend the book by Cannas da Silva \cite{can}. The rich subject of $T^n$ actions on symplectic manifolds is discussed at length by Guillemin~\cite{gui}.

\subsection{Polynomiality of Weil--Petersson volumes} \label{polynomiality}

We now consider Mirzakhani's construction of a setup in which Corollary~\ref{volumevar} may be used to produce Weil--Petersson volumes~\cite{mir1}. This allows us to prove that $V_{g,n}(\LL)$ is a polynomial and, furthermore, that its coefficients store intersection numbers on the moduli space of curves $\M_{g,n}$. We start by considering the space
\[
\widehat{{\mathcal M}}_{g,n} = \left\{ (X, p_1, p_2, \ldots, p_n) \; \middle \vert \; \begin{array}{l} X \text{ is a genus } g \text{ hyperbolic surface with} n \text{ geodesic boundary} \\ \text{components } \beta_1, \beta_2, \ldots, \beta_n \text{ and } p_k \in \beta_k \text{ for all } k \end{array} \right\}.
\]
There is a $T^n$ action on this space, where the $k$th copy of $S^1$ moves the point $p_k$ along the boundary $\beta_k$ at a constant speed in the direction given by the orientation of the surface.

We now show that $\widehat{\mathcal M}_{g,n}$ has a $T^n$ invariant symplectic structure. Fix a tuple $\gamma = (\gamma_1, \gamma_2, \ldots, \gamma_n)$ of homotopy classes of disjoint simple closed curves on the surface $\Sigma_{g,2n}$ with genus $g$ and $2n$ labelled boundary components such that $\gamma_k$ bounds a pair of pants with the boundaries labelled $2k-1$ and $2k$. Since mapping classes act on homotopy classes of curves, elements of $\textup{Mod}_{g,2n}$ act on $\gamma$ componentwise. Now define
\[
{\mathcal M}_{g,2n}^* = \{ (X, \eta_1, \eta_2, \ldots, \eta_n) \mid X \in {\mathcal M}_{g,2n}(\mathbf{0}) \text{ and } (\eta_1, \eta_2, \ldots, \eta_n) \in \textup{Mod}_{g,2n} \cdot \gamma\}.
\]
Equivalently, we can use the definition ${\mathcal M}_{g,2n}^* = {\mathcal T}_{g,2n}(\mathbf{0}) / \textup{Stab}(\gamma)$, where the stabiliser
\[
\textup{Stab}(\gamma) = \{[\phi] \in \textup{Mod}_{g,2n} \mid \phi(\gamma_k) \text{ is homotopic to } \gamma_k \text{ for all } k\} \leq \textup{Mod}_{g,2n}
\]
acts on the Teichm\"{u}ller space in the usual way. Since the Weil--Petersson symplectic form on the Teichm\"{u}ller space ${\mathcal T}_{g,2n}(\mathbf{0})$ is invariant under the action of the mapping class group, it must also be invariant under $\textup{Stab}(\gamma)$. Therefore, it descends to a symplectic form on ${\mathcal M}_{g,2n}^*$.

There is a natural map $f: \widehat{{\mathcal M}}_{g,n} \to {\mathcal M}_{g,2n}^*$ which is easy to describe. Simply take $(X, p_1, p_2, \ldots, p_n)$ where $X \in {\mathcal M}_{g,n}(\LL)$ and, to the $k$th boundary component, glue in a pair of pants with two cusps labelled $2k-1$ and $2k$ and a boundary component of length $L_k$. Of course, this can be done in infinitely many ways and we choose the unique way such that the seam from the cusp labelled $2k$ meets the point $p_k$. The map $f$ can be used to pull back the symplectic form from ${\mathcal M}_{g,2n}^*$ to $\widehat{{\mathcal M}}_{g,n}$, where it is invariant under the $T^n$ action. Furthermore, by the definition of the Weil--Petersson symplectic form, the canonical map $\ell^{-1}(\LL) / T^n \to {\mathcal M}_{g,n}(\LL)$ is a symplectomorphism, where $\ell: \widehat{\mathcal M}_{g,n} \to \mathbb{R}^n$ sends a hyperbolic surface to its boundary lengths. One may check that the $T^n$ action is the Hamiltonian flow for the moment map $\mu: \widehat{{\mathcal M}}_{g,n} \to \mathbb{R}^n$ defined by $\mu(X, p_1, p_2, \ldots, p_n) =(\frac{1}{2} L_1^2, \frac{1}{2} L_2^2, \ldots, \frac{1}{2} L_n^2)$, where $L_k$ denotes the length of the geodesic boundary component $\beta_k$.

By construction, the symplectic quotient $\mu^{-1}(\frac{1}{2} L_1^2, \frac{1}{2} L_2^2, \ldots, \frac{1}{2} L_n^2) / T^n$ is the moduli space ${\mathcal M}_{g,n}(\LL)$. As usual, the $T^n$ action gives rise to $n$ circle bundles on the symplectic quotient. Although the moment map is only regular away from $\mathbf{0}$, one obtains circle bundles ${\mathcal S}_1, {\mathcal S}_2, \ldots, {\mathcal S}_n$ on ${\mathcal M}_{g,n}(\mathbf{0})$ by taking the limit as $\LL \to \mathbf{0}$. Mirzakhani proved the following fact concerning the Chern classes of these circle bundles.

\begin{proposition} \label{circlebundle}
For $k = 1, 2, \ldots, n$, $c_1({\mathcal S}_k) = \psi_k \in H^2(\M_{g,n}; \mathbb{Q})$.
\end{proposition}

This proposition states that $c_1({\mathcal S}_k)$ is an element of $H^2(\M_{g,n}; \mathbb{Q})$, even though it is apparent that ${\mathcal S}_k$ is a circle bundle over the uncompactified space ${\mathcal M}_{g,n}(\mathbf{0})$. However, with a little more care, all of the previous discussion generalises to the Deligne--Mumford compactifications of the moduli spaces involved. We are now ready to state and prove one of the most important results underlying this article.

\begin{theorem}[Mirzakhani's theorem] \label{mirzakhani}
The Weil--Petersson volume $V_{g,n}(\LL)$ is given by the formula
\[
\sum_{|\bm{\alpha}| + m = 3g-3+n} \frac{(2\pi^2)^m \int_{\M_{g,n}} \psi_1^{\alpha_1} \psi_2^{\alpha_2} \cdots \psi_n^{\alpha_n} \kappa_1^{m}}{2^{|\bm{\alpha}|} \alpha_1! \alpha_2! \cdots \alpha_n! m!} L_1^{2\alpha_1} L_2^{2\alpha_2} \cdots L_n^{2\alpha_n}.
\]
\end{theorem}

\begin{proof}
We simply apply Corollary~\ref{volumevar} to the symplectic manifold $\widehat{{\mathcal M}}_{g,n}$ with the moment map $\mu$ defined above. This implies that the Weil--Petersson volume of ${\mathcal M}_{g,n}(\LL)$ for $\LL \neq \mathbf{0}$ is a polynomial in $\frac{1}{2}L_1^2, \frac{1}{2}L_2^2, \ldots, \frac{1}{2}L_n^2$. The coefficients are given by integrating products of Chern classes of certain circle bundles alongside powers of the reduced symplectic form. In the $\LL \to \mathbf{0}$ limit, Proposition~\ref{circlebundle} states that these Chern classes are precisely the psi-classes on $\M_{g,n}$. Furthermore, the reduced symplectic form converges to the usual Weil--Petersson symplectic form in the limit. All that is required now is to invoke Corollary~\ref{volumevar} and substitute $\omega = 2\pi^2 \kappa_1$, which is true in cohomology by Theorem~\ref{kappa1}.\footnote{The literature on symplectic reduction generally does not discuss the case of symplectic orbifolds. However, one can get around such problems by lifting to a manifold cover, which is possible by Theorem~\ref{manifold-cover}. This takes a little extra care, but essentially causes no problems.}
\end{proof}

We remark that Theorem~\ref{formvar} applied to this setup yields a generalisation of Theorem~\ref{kappa1}. The generalisation states that the de Rham cohomology class of the Weil--Petersson symplectic form on $\overline{\mathcal M}_{g,n}(\LL)$ satisfies 
\[
[\omega] = 2\pi^2 \kappa_1 + \frac{1}{2} L_1^2 \psi_1 + \frac{1}{2} L_2^2 \psi_2 + \cdots + \frac{1}{2} L_n^2 \psi_n \in H^2(\M_{g,n}; \mathbb{R}).
\]

Mirzakhani's theorem shows that the Weil--Petersson volume $V_{g,n}(\LL)$ is a polynomial whose coefficients store intersection numbers on $\M_{g,n}$. One of the consequences is that any meaningful statement about the volume $V_{g,n}(\LL)$ yields a meaningful statement about the intersection theory on $\M_{g,n}$, and vice versa. In this section, we have only outlined the proof of Mirzakhani's theorem, whereas the technical details may be found in Mirzakhani's original paper~\cite{mir1}.

\section{A recursion for Weil--Petersson volumes} \label{recursive-formula}

\subsection{The volume of ${\mathcal M}_{1,1}(0)$} \label{volume-M11}

One of the main obstacles in the calculation of Weil--Petersson volumes is the fact that the Fenchel--Nielsen coordinates for Teichm\"{u}ller space are not well-behaved under the action of the mapping class group. In particular, there is no concrete description for a fundamental domain of ${\mathcal M}_{g,n}(\LL)$ in ${\mathcal T}_{g,n}(\LL)$ in the general case. The workaround successfully applied by Mirzakhani~\cite{mir2} is to {\em unfold the integral} in the following way. Let $\pi: X_1 \to X_2$ be a covering map, $dv_2$ a volume form on $X_2$, and $dv_1 = \pi^*dv_2$ the pullback volume form on $X_1$. If $\pi$ is a finite covering, then for any function $f: X_1 \to \mathbb{R}$, one can construct the pushforward function $\pi_*f: X_2 \to \mathbb{R}$ defined by
\[
(\pi_*f)(y) = \sum_{x \in \pi^{-1}(y)} f(x).
\]
In fact, even if $\pi$ is an infinite covering, then the pushforward function may still exist, provided $f$ is sufficiently well-behaved. The main reason for considering this setup is the fact that, under mild integrability assumptions, we have
\[
\int_{X_1} f~dv_1 = \int_{X_2} (\pi_*f)~dv_2.
\]

We will use this strategy to calculate the volume of ${\mathcal M}_{1,1}(0)$, which will serve as a motivating example for the general case of ${\mathcal M}_{g,n}(\LL)$. For this, set $X_2 = {\mathcal M}_{1,1}(0)$ and
\[
X_1 = {\mathcal M}_{1,1}^*(0) = \{(X, \gamma) \mid X \in {\mathcal M}_{1,1}(0) \text{ and $\gamma$ a simple closed geodesic on } X\}.
\]
Equivalently, we can use the definition ${\mathcal M}_{1,1}^*(0) = {\mathcal T}_{1,1}(0) / \textup{Stab}(\alpha)$, where $\alpha$ is a simple closed curve on the interior of the once-punctured torus. The stabiliser
\[
\textup{Stab}(\alpha) = \{[\phi] \in \textup{Mod}_{1,1} \mid \phi(\alpha) \text{ is homotopic to } \alpha\} \leq \textup{Mod}_{1,1}
\]
acts on the Teichm\"{u}ller space in the usual way. Using Fenchel--Nielsen coordinates, each $(X, \gamma) \in {\mathcal M}_{1,1}^*(0)$ can be described by the pair $(\ell, \tau)$, where $\ell$ denotes the length of $\gamma$ and $\tau$ the corresponding twist parameter. The only redundancy in this description comes from the fact that the pair $(\ell, \tau+\ell)$ may also be used to describe the same point in ${\mathcal M}_{1,1}^*(0)$. Hence, we can write
\[
{\mathcal M}_{1,1}^*(0) \cong \{ (\ell, \tau) \mid \ell \in \mathbb{R}_+ \text{ and } 0 \leq \tau \leq \ell \} / \sim,
\]
where $(\ell,0) \sim (\ell,\ell)$ for all $\ell \in \mathbb{R}_+$.

The map $\pi: {\mathcal M}_{1,1}^*(0) \to {\mathcal M}_{1,1}(0)$ is the obvious projection map defined by $\pi(X, \gamma) = X$. Through the tower of coverings ${\mathcal T}_{1,1}(0) \to {\mathcal M}_{1,1}^*(0) \to {\mathcal M}_{1,1}(0)$, we see that the Weil--Petersson form pulls back to $\pi^*\omega = d\ell \wedge d\tau$ on the intermediate cover ${\mathcal M}_{1,1}^*(0)$. Let $\ell: {\mathcal M}_{1,1}^*(0) \to \mathbb{R}$ be the geodesic length function so that $\ell(X, \gamma)$ equals the length of $\gamma$ on $X$. Unfolding the integral and using the description above for ${\mathcal M}_{1,1}^*(0)$ yields the following equalities.
\[
\int_{{\mathcal M}_{1,1}(0)} \sum_{\pi(Y) = X} f(\ell(Y))\,dX = \int_{{\mathcal M}_{1,1}^*(0)} f(\ell(Y))\,dY = \int_0^\infty \!\! \int_0^\ell f(\ell)\,d\tau\,d\ell
\]

Therefore, in order to obtain the volume of ${\mathcal M}_{1,1}(0)$, we would like an identity of the form
\[
\sum_{\pi(Y) = X} f(\ell(Y)) = 1,
\]
valid for all $X \in {\mathcal M}_{1,1}(0)$. Note that the summation is over the set of simple closed geodesics on $X$. Such an identity had been discovered by McShane~\cite{mcs1} prior to Mirzakhani's work on Weil--Petersson volumes.

\begin{theorem}[McShane identity] \label{mcshane}
If $X$ is a hyperbolic torus with one cusp, then
\[
\sum_{\gamma} \frac{2}{1+\exp \ell(\gamma)} = 1.
\]
Here, the summation is over the set of simple closed geodesics on $X$ and $\ell(\gamma)$ denotes the length of $\gamma$.
\end{theorem}

So to complete the calculation of the volume of ${\mathcal M}_{1,1}(0)$, take $f(\ell) = \frac{2}{1+\exp \ell}$.
\[
\int_{{\mathcal M}_{1,1}(0)} 1\,dX = \int_0^\infty \!\! \int_0^\ell \frac{2}{1+\exp \ell}\,d\tau\,d\ell = \int_0^\infty \frac{2\ell}{1+\exp \ell}\,d\ell = \frac{\pi^2}{6}
\]
However, recall that the case of ${\mathcal M}_{1,1}(0)$ is exceptional in the sense that a generic point of the moduli space is an orbifold point with orbifold group $\mathbb{Z} / 2\mathbb{Z}$. This fact arises from the existence of the elliptic involution on every hyperbolic torus with one cusp. Since we consider orbifold volumes in this article, it is necessary to divide the integral calculation above by two. Therefore, we finally have the result
\[
V_{1,1}(0) = \frac{\pi^2}{12}.
\]

\subsection{McShane identities}

In order to unfold the integral required to calculate $V_{g,n}(\LL)$, it is necessary to obtain a more general version of McShane's identity. The following generalisation is due to Mirzakhani~\cite{mir2}.

\begin{theorem}[Generalised McShane identity] \label{generalised-mcshane}
On a hyperbolic surface with geodesic boundary components $\beta_1, \beta_2, \ldots, \beta_n$ of lengths $L_1, L_2, \ldots, L_n$, respectively,
\[
\sum_{(\alpha_1, \alpha_2)} D(L_1, \ell(\alpha_1), \ell(\alpha_2)) + \sum_{k=2}^n \sum_{\gamma} R(L_1, L_k, \ell(\gamma)) = L_1.
\]
Here, the first summation is over unordered pairs $(\alpha_1, \alpha_2)$ of simple closed geodesics which bound a pair of pants with $\beta_1$, while the second summation is over simple closed geodesics $\gamma$ which bound a pair of pants with $\beta_1$ and $\beta_k$. The functions $D: \mathbb{R}^3 \to \mathbb{R}$ and $R: \mathbb{R}^3 \to \mathbb{R}$ are given by the equations
\[
D(x, y, z) = 2 \log \left( \frac{e^{\frac{x}{2}} + e^{\frac{y+z}{2}}}{e^{-\frac{x}{2}} + e^{\frac{y+z}{2}}} \right) \text{ and } R(x, y, z) = x - \log  \left( \frac{\cosh \frac{y}{2} + \cosh \frac{x+z}{2}}{\cosh \frac{y}{2} + \cosh \frac{x-z}{2}} \right).
\]
\end{theorem}

The main idea behind the proof is to consider, for each point $x \in \beta_1$, the geodesic $\gamma_x$ which meets $\beta_1$ orthogonally at $x$. If we start at $x$ and walk along $\gamma_x$, then one of the following situations must arise.
\begin{enumerate} [topsep=0pt,itemsep=1pt]
\item The geodesic $\gamma_x$ intersects itself. \label{case1}
\item The geodesic $\gamma_x$ intersects $\beta_1$ without intersecting itself. \label{case2}
\item The geodesic $\gamma_x$ intersects $\beta_k$ for $2 \leq k \leq n$ without intersecting itself. \label{case3}
\item The geodesic $\gamma_x$ never intersects itself or a boundary component. \label{case4}
\end{enumerate}
We now use this observation to construct a map from a subset $\beta_1^* \subseteq \beta_1$ to the set
\[
{\mathcal P}_1 = \left\{ \begin{array}{c} \text{embedded hyperbolic pairs of pants, one of} \\ \text{whose geodesic boundary components is } \beta_1 \end{array} \right\}.
\]
Note that the generalised McShane identity is not a summation over simple closed geodesics, but over ${\mathcal P}_1$. It just so happens that on a once-punctured torus, the two notions coincide. In fact, McShane's original identity --- see Theorem~\ref{mcshane} --- may be recovered by using $(g,n) = (1,1)$, dividing both sides of the identity by $L_1$, and taking the $L_1 \to 0$ limit.

In cases (\ref{case1}) and (\ref{case2}), consider the union of $\beta_1$ and the geodesic path $\gamma_x$ from $x$ to the intersection point. For $\epsilon > 0$ sufficiently small, the $\epsilon$-neighbourhood of this embedded graph is topologically a pair of pants. By taking geodesic representatives in the homotopy classes of the boundary components, we obtain an embedded hyperbolic pair of pants, one of whose geodesic boundary components is $\beta_1$. Let $f(x) \in {\mathcal P}_1$ denote this pair of pants.

In case (\ref{case3}), consider the union of $\beta_1$, $\beta_k$ and the geodesic path $\gamma_x$ from $x$ to the intersection point. Again, for $\epsilon > 0$ sufficiently small, the $\epsilon$-neighbourhood of this embedded graph is topologically a pair of pants. By taking geodesic representatives in the homotopy classes of the boundary components, we obtain an embedded hyperbolic pair of pants, one of whose geodesic boundary components is $\beta_1$. Let $f(x) \in {\mathcal P}_1$ denote this pair of pants.

Thus, we have defined a function $f: \beta_1^* \to {\mathcal P}_1$, where $\beta_1^*$ is the set consisting of those points in $\beta_1$ for which cases (\ref{case1}), (\ref{case2}) or (\ref{case3}) occur. The points in $\beta_1$ for which $f$ is undefined are those for which case (\ref{case4}) occurs. A result due to Birman and Series~\cite{bir-ser} states that the union of all complete simple geodesics on a closed hyperbolic surface has Hausdorff dimension one. By doubling the surface along its boundary, we can generalise the statement to hyperbolic surfaces with boundary and complete simple geodesics perpendicular to the boundary. Hence, we may deduce that $\beta \setminus \beta_1^*$ is a subset of zero measure with respect to the hyperbolic line element $\mu$ on $\beta_1$. In fact, Mirzakhani~\cite{mir2} shows that it is homeomorphic to the union of a Cantor set and countably many isolated points. The upshot of this discussion is the equation
\[
\sum_{P \in {\mathcal P}_1} \mu(f^{-1}(P)) = L_1.
\]

Theorem~\ref{generalised-mcshane} now follows from a couple of simple facts.

\begin{lemma} \label{D-and-R} ~
\begin{itemize} [topsep=0pt,itemsep=1pt]
\item If $P \in {\mathcal P}_1$ is bound by $\beta_1$ and two simple closed geodesics $\alpha_1$ and $\alpha_2$, then
\[
\mu(f^{-1}(P)) = D(L_1, \ell(\alpha_1), \ell(\alpha_2)).
\]
\item If $P \in {\mathcal P}_1$ is bound by $\beta_1$, $\beta_k$ and a simple closed geodesic $\gamma$, then
\[
\mu(f^{-1}(P)) = R(L_1, L_k, \ell(\gamma)).
\]
\end{itemize}
\end{lemma}

Note that the calculation of $\mu(f^{-1}(P))$ is local in the sense that it depends only on the geometry of $P$ and not on the geometry of the entire surface. To a hyperbolic pair of pants with geodesic boundary components $\alpha, \beta, \gamma$, we associate four distinguished points on $\alpha$. Such a pair of pants necessarily contains exactly four complete simple geodesics which meet $\alpha$ orthogonally exactly once and are disjoint from $\beta$ and $\gamma$. One of them intersects $\alpha$ at $B_1$ and spirals around $\beta$ one way, while another intersects $\alpha$ at $B_2$ and spirals around $\beta$ the other way. Similarly, one of them intersects $\alpha$ at $C_1$ and spirals around $\gamma$ one way, while another intersects $\alpha$ at $C_2$ and spirals around $\gamma$ the other way. Note that the orientation reversing isometry which reflects the pair of pants through its seams interchanges $B_1$ with $B_2$ and $C_1$ with $C_2$. This is encapsulated in the schematic diagram below.

\vspace{6pt}
\begin{center}
\begin{overpic}[unit=1mm]{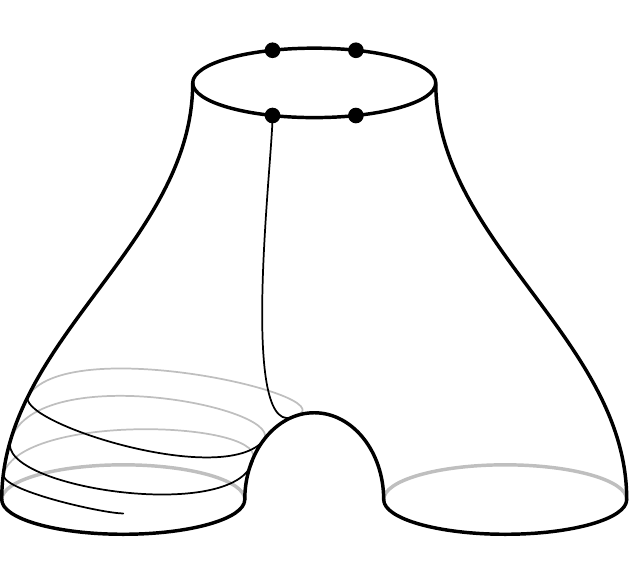}
\put(49,78.5){$\alpha$}
\put(17,2){$\beta$}
\put(81,2){$\gamma$}
\put(36,69){$B_1$}
\put(41,88){$B_2$}
\put(55,69){$C_1$}
\put(55,88){$C_2$}
\end{overpic}
\end{center}

Suppose that $x \in \alpha$ lies on the interval $B_1C_1$ which does not include $B_2$ and $C_2$ or on the interval $B_2C_2$ which does not include $B_1$ and $C_1$. Then the geodesic $\gamma_x$ will intersect itself or $\alpha$, so case (\ref{case1}) or (\ref{case2}) occurs. Therefore, we define $D(\ell(\alpha), \ell(\beta), \ell(\gamma))$ to be twice the length of the interval $B_1C_1$ which does not include $B_2$ and $C_2$. Now suppose that $x \in \alpha$ lies on the interval $B_1B_2$ which does not include $C_1$ and $C_2$. Then the geodesic $\gamma_x$ will intersect $\beta$, so case (\ref{case3}) occurs. Therefore, we define $R(\ell(\alpha), \ell(\beta), \ell(\gamma))$ to be the length of the interval $C_1C_2$ which includes $B_1$ and $B_2$. The proof of Lemma~\ref{D-and-R} now follows from these definitions.

All that remains is to explicitly compute the functions $D$ and $R$. In the universal cover, the value of $D(\ell(\alpha), \ell(\beta), \ell(\gamma))$ is twice the distance between the projection of $\beta$ and $\gamma$ on $\alpha$ and the value of $R(\ell(\alpha), \ell(\beta), \ell(\gamma))$ is $\ell(\alpha)$ minus the length of the projection of $\gamma$ on $\alpha$. We do not complete the calculation here but remark that it can be carried out by applying some elementary results from hyperbolic trigonometry.\footnote{For a valuable reference on hyperbolic trigonometry and the hyperbolic geometry of surfaces, see Buser's {\em Geometry and spectra of compact Riemann surfaces}~\cite{bus}.}

There are now many variations on the McShane theme.
\begin{itemize} [topsep=0pt,itemsep=1pt]
\item Bowditch used the notion of Markoff triples to give an alternative proof of McShane's identity for the once-punctured torus~\cite{bow}.

\item Akiyoshi, Miyachi and Sakuma produced variants of McShane's identity for quasi-fuchsian punctured surface groups and hyperbolic punctured surface bundles over the circle~\cite{aki-miy-sak}.

\item McShane determined an identity for simple geodesics on a closed hyperbolic surface of genus two~\cite{mcs2}. This work capitalises on the existence of the hyperelliptic involution and Weierstrass points on a genus two surface.

\item Tan, Wong and Zhang gave a generalisation of McShane's identity to hyperbolic cone surfaces, where all cone points have angles bounded above by $\pi$~\cite{tan-won-zha1}. They also found variations concerning representations of punctured torus groups to $SL(2, \mathbb{C})$~\cite{tan-won-zha2} and also classical Schottky groups~\cite{tan-won-zha3}.

\item Luo and Tan have recently found a McShane identity for all closed hyperbolic surfaces~\cite{luo-tan}. Their proof draws some inspiration from Calegari's elegant and unified treatment of the following two results --- namely, the identities of Basmajian~\cite{bas} and Bridgeman~\cite{bri-kah}. These share a similar flavour with McShane identities, but pertain to hyperbolic manifolds of arbitrary dimension.
\end{itemize}

\begin{theorem}
Let $M$ be a compact hyperbolic $n$-manifold with totally geodesic boundary $\partial M$ and let $(\ell_1, \ell_2, \ell_3, \ldots)$ denote the lengths of the orthogeodesics of $M$, listed with multiplicity.\footnote{An {\em orthogeodesic} of $M$ is a geodesic arc which is perpendicular to $\partial M$ at its endpoints.} There exist functions $A_n$ and $V_n$ depending only on $n$ such that the following identities hold~\cite{cal}.
\[
\textup{area } \partial M = \sum A_n(\ell_i) \quad \text{and} \quad \textup{volume } M= \sum V_n(\ell_i)
\]
\end{theorem}

\subsection{Mirzakhani's recursion}

In this section, we prove the following formula for Weil--Petersson volumes, originally due to Mirzakhani~\cite{mir2}. We use the convention that $V_{0,1}(L_1) = 0$, $V_{0,2}(L_1, L_2) = 0$ and $V_{0,3}(L_1, L_2, L_3) = 1$.

\begin{theorem}[Mirzakhani's recursion] \label{mirzakhani-recursion} The Weil--Petersson volumes satisfy the following equation for $2g+n > 3$.
\begin{align*}
2 \frac{\partial}{\partial L_1} L_1 V_{g,n}(\LL) = & \int_0^\infty \!\! \int_0^\infty xy\,H(x+y, L_1)\,V_{g-1,n+1}(x, y, \widehat{\LL}) \,dx\,dy \\
& + \sum_{\substack{g_1+g_2=g \\ I \sqcup J = [2, n]}} \int_0^\infty \!\! \int_0^\infty xy\,H(x+y, L_1)\,V_{g_1,|I|+1}(x, \LL_I) \,V_{g_2,|J|+1}(y, \LL_J) \,dx\,dy \\
& + \sum_{k=2}^n \int_0^\infty x \, [H(x, L_1+L_k) + H(x, L_1-L_k)]\, V_{g,n-1}(x, \widehat{\LL}_k) \,dx
\end{align*}
Here, we have used the notation $\widehat{\LL} = (L_2, L_3, \ldots, L_n)$, $\LL_I = (L_{i_1}, L_{i_2}, \ldots, L_{i_m})$ for $I = \{i_1, i_2, \ldots, i_m\}$, and $\widehat{\LL}_k = (L_2, \ldots, \widehat{L}_k, \ldots, L_n)$ where the hat denotes omission. The function $H: \mathbb{R}^2 \to \mathbb{R}$ is defined by
\[
H(x, y) = \frac{1}{1 + \exp \frac{x+y}{2}} + \frac{1}{1 + \exp \frac{x-y}{2}}.
\]
\end{theorem}

The proof uses the calculation of $V_{1,1}(0)$ from Section~\ref{volume-M11} as a model. Our point of departure is the generalised McShane identity --- see Theorem~\ref{generalised-mcshane} --- which we rewrite in the following way.
\[
{\mathcal D}_{\textup{con}}(X) + \sum_{\substack{g_1+g_2=g \\ I \sqcup J = [2, n]}} {\mathcal D}_{g_1, I}(X) + \sum_{k=2}^n {\mathcal R}_k(X) = L_1
\]
We have grouped the left hand side into terms of three distinct types.
\begin{itemize} [topsep=0pt,itemsep=1pt]
\item The first type is
\[
{\mathcal D}_{\textup{con}}(X) = \sum_{(\alpha_1, \alpha_2)} D(L_1, \ell(\alpha_1), \ell(\alpha_2)),
\]
where the summation is over unordered pairs $(\alpha_1, \alpha_2)$ of simple closed geodesics which bound a pair of pants with $\beta_1$, whose complement is a connected surface.

\item The second type is
\[
{\mathcal D}_{g_1, I}(X) = \sum_{(\alpha_1, \alpha_2)} D(L_1, \ell(\alpha_1), \ell(\alpha_2)),
\]
where the summation is over unordered pairs $(\alpha_1, \alpha_2)$ of simple closed geodesics which bound a pair of pants with $\beta_1$, whose complement is a disconnected surface. We require that one component of this disconnected surface has genus $g_1$ and includes only those boundary components from the original surface labelled by elements of $I$.

\item The third type is
\[
{\mathcal R}_k(X) = \sum_{\gamma} R(L_1, L_k, \ell(\gamma)),
\]
where the summation is over simple closed geodesics $\gamma$ which bound a pair of pants with $\beta_1$ and $\beta_k$.
\end{itemize}

The rationale for expressing the generalised McShane identity in this way is that each term is now a summation over a mapping class group orbit. This is due to the fact that two sets of disjoint simple closed curves on a surface are in the same mapping class group orbit if and only if their complements have the same topological type and labelling of boundary components. Sums over mapping class group orbits can be expressed as pushforwards of functions on appropriate covers of the moduli space. And these are the functions which we are able to integrate over the moduli space itself.

Now take the generalised McShane identity and integrate both sides over the moduli space ${\mathcal M}_{g,n}(\LL)$.
\[
L_1 V_{g,n}(\LL) = \int_{{\mathcal M}_{g,n}(\LL)} {\mathcal D}_{\textup{con}}(X)\,dX + \sum_{\substack{g_1+g_2=g \\ I \sqcup J = [2, n]}} \int_{{\mathcal M}_{g,n}(\LL)} {\mathcal D}_{g_1,I}(X)\,dX + \sum_{k=2}^n \int_{{\mathcal M}_{g,n}(\LL)} {\mathcal R}_k(X)\,dX
\]
From the previous discussion, we know that it is possible to unfold each of the integrals using the strategy employed in Section~\ref{volume-M11}. For example, let us concentrate on the term
\[
\int_{{\mathcal M}_{g,n}(\LL)} {\mathcal R}_k(X)\,dX.
\]

In order to unfold the integral, recall that we require a covering map $\pi: X_1 \to X_2$, a volume form $dv_2$ on $X_2$, and the pullback volume form $dv_1 = \pi^*dv_2$ on $X_1$. For this, set $X_2 = {\mathcal M}_{g,n}(\LL)$ and
\[
X_1 = {\mathcal M}_{g,n}^*(\LL) = \left\{ (X, \gamma) \; \middle \vert \; \begin{array}{l} X \in M_{g,n}(\LL) \text{ and } \gamma \text{ a simple closed geodesic on } \\ X \text{ which bounds a pair of pants with } \beta_1 \text{ and } \beta_k \end{array} \right\}.
\]
Equivalently, we can use the definition ${\mathcal M}_{g,n}^*(\LL) = {\mathcal T}_{g,n}(\LL) / \textup{Stab}(\alpha)$, where $\alpha$ is a simple closed curve on the surface $\Sigma_{g,n}$ which bounds a pair of pants with the boundary components labelled 1 and $k$. The stabiliser
\[
\textup{Stab}(\alpha) = \{ [\phi] \in \textup{Mod}_{g,n} \mid \phi(\alpha) \text{ is homotopic to } \alpha \} \leq \textup{Mod}_{g,n}
\]
acts on the Teichm\"{u}ller space in the usual way. Using Fenchel--Nielsen coordinates, each $(X, \gamma) \in {\mathcal M}_{g,n}^*(\LL)$ can be described by the triple $(\ell, \tau, \widehat{X})$, where $\ell$ denotes the length of $\gamma$ and $\tau$ the corresponding twist parameter. The surface $\widehat{X} \in {\mathcal M}_{g,n-1}(\ell, \widehat{\LL}_k)$ is simply the complement of the pair of pants bound by $\beta_1$, $\beta_k$ and $\gamma$. The only redundancy in this description comes from the fact that the triple $(\ell, \tau+\ell, \widehat{X})$ may also be used to describe the same point in ${\mathcal M}_{g,n}^*(\LL)$. Hence, we can write
\[
{\mathcal M}_{g,n}^*(\LL) \cong \{ (\ell, \tau, \widehat{X}) \mid \ell \in \mathbb{R}_+, 0 \leq \tau \leq \ell \text{ and } \widehat{X} \in {\mathcal M}_{g,n-1}(\ell, \widehat{\LL}_k) \} / \sim,
\]
where $(\ell, \tau, \widehat{X}) \sim (\ell, \tau + \ell, \widehat{X})$.

The map $\pi: {\mathcal M}_{g,n}^*(\LL) \to {\mathcal M}_{g,n}(\LL)$ is the obvious projection map defined by $\pi(X, \gamma) = X$. Through the tower of coverings ${\mathcal T}_{g,n}(\LL) \to {\mathcal M}_{g,n}^*(\LL) \to {\mathcal M}_{g,n}(\LL)$, we see that the Weil--Petersson form pulls back to $\pi^*\omega = d\ell \wedge d\tau \wedge \widehat{\omega}$ on the intermediate cover ${\mathcal M}_{g,n}^*(\LL)$, where $\widehat{\omega}$ is the Weil--Petersson form on the lower dimensional moduli space ${\mathcal M}_{g,n-1}(\ell, \widehat{\LL}_k)$. Let $\ell: {\mathcal M}_{g,n}^*(\LL) \to \mathbb{R}$ be the geodesic length function so that $\ell(X, \gamma)$ equals the length of $\gamma$ on $X$. Unfolding the integral and using the description above for ${\mathcal M}_{g,n}^*(\LL)$ yields the following equalities.
\begin{align*}
\int_{{\mathcal M}_{g,n}(\LL)} {\mathcal R}_k(X)\,dX =& \int_{{\mathcal M}_{g,n}(\LL)} \sum_{\pi(Y) = X} R(L_1, L_k, \ell(Y))\,dX = \int_{{\mathcal M}_{g,n}^*(\LL)} R(L_1, L_k, \ell(Y))\,dY \\
=& \int_0^\infty \!\! \int_0^\ell \int_{{\mathcal M}_{g,n-1}(\ell, \widehat{\LL}_k)} R(L_1, L_k, \ell)\,\widehat{\omega}\,d\tau\,d\ell = \int_0^\infty x \,R(L_1, L_k, x) \,V_{g,n-1}(x, \widehat{\LL}_k)\,dx
\end{align*}

For the other terms in the generalised McShane identity, although the details may be different, the argument remains the same. After unfolding each of the integrals and summing, the end result is the following formula.
\begin{align*}
L_1 V_{g,n}(\LL) = & \frac{1}{2} \int_0^\infty \!\! \int_0^\infty xy\,D(L_1,x,y)\,V_{g-1,n+1}(x, y, \widehat{\LL}) \,dx\,dy \\
&+ \frac{1}{2} \sum_{\substack{g_1+g_2=g \\ I \sqcup J = [2, n]}} \int_0^\infty \!\! \int_0^\infty xy\,D(L_1,x,y)\,V_{g_1,|I|+1}(x, \LL_I) \,V_{g_2,|J|+1}(y, \LL_J) \,dx\,dy \\
&+ \sum_{k=2}^n \int_0^\infty x\,R(L_1,L_k,x)\, V_{g,n-1}(x, \widehat{\LL}_k) \,dx
\end{align*}
Note that the factor of $\frac{1}{2}$ in front of the first two terms of the right hand side is to account for the twofold symmetry between $x$ and $y$. Theorem~\ref{mirzakhani-recursion} expresses the recursion in a more useful form, which is obtained by applying $2\frac{\partial}{\partial L_1}$ to both sides of this equation.

\subsection{Applications of Mirzakhani's recursion} \label{applications}

The mechanism behind Mirzakhani's recursion is based on removing pairs of pants from the surface $\Sigma_{g,n}$ which contain at least one boundary component. Therefore, the calculation of any Weil--Petersson volume can be reduced to the base cases $V_{0,3}(L_1, L_2, L_3) = 1$ and $V_{1,1}(L_1) = \frac{1}{48}(L_1^2 + 4\pi^2)$. In the practical application of Mirzakhani's recursion, we require the following explicit integral calculations.
\begin{align*}
\int_0^\infty x^{2k-1} H(x, t)\,dx &= F_{2k-1}(t) \\
\int_0^\infty \!\! \int_0^\infty x^{2a-1} y^{2b-1} H(x+y,t)\,dx\,dy &= \frac{(2a-1)! (2b-1)!}{(2a+2b-1)!} F_{2a+2b-1}(t)
\end{align*}
Here, $F_{2k-1}(t)$ is the following even polynomial of degree $2k$ in $t$.where the coefficient of $t^{2m}$ is a rational multiple of $\pi^{2k-2m}$.
\[
F_{2k-1}(t) = (2k-1)! \sum_{i=0}^{k} \frac{\zeta(2i) (2^{2i+1}-4)}{(2k-2i)!} t^{2k-2i}
\]
We list below the polynomials $F_{2k-1}(t)$ for the first few values of $k$.
\begin{align*}
F_1(t) &= \frac{t^2}{2} + \frac{2\pi^2}{3} \\
F_3(t) &= \frac{t^4}{4} + 2\pi^2 t^2 + \frac{28\pi^4}{15} \\
F_5(t) &= \frac{t^6}{6} + \frac{10 \pi^2 t^4}{3} + \frac{56 \pi^4 t^2}{3} + \frac{992 \pi^6}{63} \\
F_7(t) &= \frac{t^8}{8} + \frac{14 \pi^2 t^6}{3} + \frac{196 \pi^4 t^4}{3} + \frac{992 \pi^6 t^2}{3} + \frac{4064 \pi^8}{15}
\end{align*}
As an example of Mirzakhani's recursion being used to calculate Weil--Petersson volumes, consider the following calculation of $V_{1,2}(L_1, L_2)$.

\begin{example}
For $(g,n) = (1,2)$, we obtain a contribution from only two terms on the right hand side of Mirzakhani's recursion --- one involving $V_{0,3}$ and the other involving $V_{1,1}$. This corresponds to the fact that removing a pair of pants from the surface $\Sigma_{1,2}$ which contains at least one boundary component must leave $\Sigma_{0,3}$ or $\Sigma_{1,1}$.

\begin{align*}
&\, 2 \frac{\partial}{\partial L_1} L_1 V_{1,2}(L_1, L_2) \\
=\,& \int_0^\infty \!\! \int_0^\infty xy\,H(x+y,L_1)\,V_{0,3}(x, y, L_2)\,dx\,dy + \int_0^\infty x [H(x,L_1+L_2) + H(x,L_1-L_2)]\,V_{1,1}(x)\,dx \\
=\,& \int_0^\infty \!\! \int_0^\infty xy\,H(x+y,L_1)\,dx\,dy + \int_0^\infty x [H(x,L_1+L_2) + H(x,L_1-L_2)] \left(\frac{x^2 + 4\pi^2}{48}\right) \,dx \\
=\,& \frac{F_3(L_1)}{6} + \frac{F_3(L_1+L_2) + F_3(L_1-L_2)}{48} + \frac{\pi^2 F_1(L_1+L_2) + \pi^2 F_1(L_1-L_2)}{12} \\
=\,& \frac{5L_1^4}{96} + \frac{L_1^2L_2^2}{16} + \frac{L_2^4}{96} + \frac{\pi^2L_1^2}{2} + \frac{\pi^2L_2^2}{6} + \frac{\pi^4}{2}
\end{align*}

Now integrate with respect to $L_1$ and divide by $2L_1$ to obtain the desired result. Observe that no constant of integration appears, since $V_{1,2}(0, 0)$ is finite, as noted in Section~\ref{early-results}.
\[
V_{1,2}(L_1, L_2) = \frac{L_1^4}{192} + \frac{L_1^2L_2^2}{96} + \frac{L_2^4}{192} + \frac{\pi^2 L_1^2}{12} + \frac{\pi^2L_2^2}{12} + \frac{\pi^4}{4}
\]
\end{example}

Mirzakhani's recursion can be used to provide an alternative proof of the following fact, which is a direct corollary of Theorem~\ref{mirzakhani}.

\begin{corollary} \label{volume-polynomial}
The Weil--Petersson volume $V_{g,n}(\LL)$ is an even symmetric polynomial in $L_1, L_2, \ldots, L_n$ of degree $6g-6+2n$. Furthermore, the coefficient of $L_1^{2\alpha_1} L_2^{2\alpha_2} \cdots L_n^{2\alpha_n}$ is a rational multiple of $\pi^{6g-6+2n-2|\bm{\alpha}|}$.
\end{corollary}

The symmetry of $V_{g,n}(\LL)$ is a consequence of the symmetry of the boundary labels. However, the symmetry is not present in Mirzakhani's recursion, which treats one of the boundary components as distinguished. The remainder of Corollary~\ref{volume-polynomial} can be proven with a straightforward application of induction on the value of $2g-2+n$.

Mirzakhani's theorem --- see Theorem~\ref{mirzakhani} --- shows that $V_{g,n}(\LL)$ is a polynomial whose coefficients store information about the intersection theory on $\M_{g,n}$. In fact, all psi-class intersection numbers on $\M_{g,n}$ can be recovered from the top degree part of $V_{g,n}(\LL)$ alone. On the other hand, Mirzakhani's recursion --- see Theorem~\ref{mirzakhani-recursion} --- shows that the Weil--Petersson volume $V_{g,n}(\LL)$ can be calculated in an explicit manner. So the conjunction of these two results provides an algorithm to compute all psi-class intersection numbers on $\M_{g,n}$. Thus, Mirzakhani was able to give a new proof of the Witten--Kontsevich theorem~\cite{mir1}. Although several proofs of the Witten--Kontsevich theorem now exist, there are three novel features of Mirzakhani's proof. First, she proved it by directly verifying the Virasoro constraints. Second, her proof was the first to appear which did not make explicit use of a matrix model. Third, her work uses hyperbolic geometry in a fundamental way.

Further mileage can be obtained from integration over moduli spaces of hyperbolic surfaces. For example, Mirzakhani has applied this technique to obtain the following result concerning the number of simple closed geodesics of bounded length on a hyperbolic surface~\cite{mir3}.

\begin{theorem}
For $\gamma$ a simple closed curve on $X \in {\mathcal M}_{g,n}(\mathbf{0})$, let $s(X, \gamma, N)$ denote the number of simple closed geodesics in the mapping class group orbit of $\gamma$ whose length is at most $N$. Then
\[
\lim_{N \to \infty} \frac{s(X, \gamma, N)}{N^{6g-6+2n}} = \frac{c(\gamma) B(X)}{\int_{{\mathcal M}_{g,n}(\mathbf{0})} B(X)}.
\]
Here, $c(\gamma) \in \mathbb{Q}$ depends only on the topological type of $\gamma$ and $B(X)$ is the volume of the unit ball centred at $X$ in the space of measured geodesic laminations.
\end{theorem}

\section{Limits of Weil--Petersson volumes} \label{WP-limits}

\subsection{Hyperbolic cone surfaces} \label{cone-surfaces}

Mirzakhani's recursion --- see Theorem~\ref{mirzakhani-recursion} --- allows us to calculate Weil--Petersson volumes explicitly. The table in Appendix~\ref{WP-table} contains $V_{g,n}(\LL)$ for various values of $g$ and $n$. These data suggest the striking observation that $V_{g,1}(2\pi i) = 0$. This statement does indeed hold true for all positive integers $g$, which indicates that the Weil--Petersson volume polynomials display interesting behaviour when the lengths are formally set to $2\pi i$. Further investigation yields the following results~\cite{do-nor}.

\begin{theorem}[String and dilaton equations for Weil--Petersson volumes] \label{string-dilaton-WP}
For $2g-2+n > 0$, the Weil--Petersson volumes satisfy the following relations.
\begin{align*}
V_{g,n+1}(\LL, 2\pi i) &= \sum_{k=1}^n \int_0^{L_k} L_k \,V_{g,n}(\LL)\,dL_k \\
\frac{\partial V_{g,n+1}}{\partial L_{n+1}}(\LL, 2\pi i) &= 2\pi i \,(2g-2+n) \,V_{g,n}(\LL)
\end{align*}
\end{theorem}

These equations must follow from Mirzakhani's recursion since it uniquely determines all Weil--Petersson volumes. The proof based on this observation is rather unwieldy and not so transparent~\cite{do1}. An alternative proof expresses the string and dilaton equations as relations between the coefficients of $V_{g,n+1}(\LL, L_{n+1})$ and of $V_{g,n}(\LL)$. By Mirzakhani's theorem --- see Theorem~\ref{mirzakhani} --- this translates to relations between intersection numbers on $\M_{g,n+1}$ and on $\M_{g,n}$. Thus, we may equivalently write Theorem~\ref{string-dilaton-WP} in the following way.
\begin{align*}
\sum_{j=0}^m (-1)^j \binom{m}{j} \int_{\M_{g,n+1}} \psi_1^{\alpha_1} \cdots \psi_n^{\alpha_n} \psi_{n+1}^j \kappa_1^{m-j} &= \sum_{k=1}^n \int_{\M_{g,n}} \psi_1^{\alpha_1} \cdots \psi_k^{\alpha_k-1} \cdots \psi_n^{\alpha_n} \kappa_1^m \\
\sum_{j=0}^m (-1)^j \binom{m}{j} \int_{\M_{g,n+1}} \psi_1^{\alpha_1} \cdots \psi_n^{\alpha_n} \psi_{n+1}^{j+1} \kappa_1^{m-j} &= (2g-2+n) \int_{\M_{g,n}} \psi_1^{\alpha_1} \cdots \psi_n^{\alpha_n} \kappa_1^m
\end{align*}
Observe that $\psi_k$ on the left hand side refers to the psi-class on $\M_{g,n+1}$ while on the right hand side it refers to the psi-class on $\M_{g,n}$. These are generalisations of the string and dilaton equations --- see Theorem~\ref{string-dilaton} --- which correspond to the case $m = 0$. They may be proven using standard arguments from algebraic geometry~\cite{do-nor}. The succinct statement of Theorem~\ref{string-dilaton-WP} indicates that the Weil--Petersson volume polynomial $V_{g,n}(\LL)$ provides a useful way to package intersection numbers on $\M_{g,n}$.

One predicts yet another approach to the string and dilaton equations, which may prove to be the most interesting. A phenomenon often occurring in hyperbolic geometry is the fact that a purely imaginary length can be interpreted as an angle. As an example, consider the work of Tan, Wong and Zhang~\cite{tan-won-zha1}, in which they show that the generalised McShane identity --- see Theorem~\ref{generalised-mcshane} --- holds for hyperbolic cone surfaces. In fact, one need only substitute $i \theta$ into the formula to represent a cone point with angle $\theta$. It follows that one can extend the definition of the Weil--Petersson volume polynomials to the case of moduli spaces of hyperbolic cone surfaces. Thus, it is tempting to think of the string and dilaton equations as describing the Weil--Petersson volume and its derivative as one of the boundary components degenerates to a cone point with angle $2\pi$ and hence, is removable. Unfortunately, a proof of the string and dilaton equations following this intuition is yet to be formalised. One of the main obstacles is the fact that the Teichm\"{u}ller theory for hyperbolic cone surfaces breaks down when cone points have angles larger than $\pi$. Indeed, on such surfaces, it ceases to be true that every homotopy class of closed curves contains a geodesic representative.

Note that Mirzakhani's recursion does not produce Weil--Petersson volumes of moduli spaces of closed hyperbolic surfaces. One application of Theorem~\ref{string-dilaton-WP} is the computation of these numbers. The following result is a direct corollary of the dilaton equation in the $n = 0$ case~\cite{do-nor}.

\begin{corollary}
The Weil--Petersson volumes of moduli spaces of closed hyperbolic surfaces satisfy the formula
\[
V_{g,0} = \frac{V_{g,1}'(2\pi i)}{2\pi i (2g-2)}.
\]
\end{corollary}

Another application of the string and dilaton equations is the computation of small genus Weil--Petersson volumes~\cite{do-nor}.

\begin{proposition} \label{genus-0-and-1}
The string equation alone uniquely determines $V_{0,n+1}(\LL, L_{n+1})$ from $V_{0,n}(\LL)$. Similarly, the string and dilaton equations together uniquely determine $V_{1,n+1}(\LL, L_{n+1})$ from $V_{1,n}(\LL)$.
\end{proposition}

The proof of Proposition~\ref{genus-0-and-1} is elementary and can be converted to algorithms for the computation of $V_{0,n}(\LL)$ and $V_{1,n}(\LL)$. These are empirically more efficient than a direct implementation of Mirzakhani's recursion, which requires $V_{0,k}$ for $3 \leq k \leq n$ for the computation of $V_{0,n+1}$.

The string and dilaton equations for Weil--Petersson volumes relate the value and derivative of $V_{g,n+1}(\LL, L_{n+1})$ evaluated at $L_{n+1} = 2\pi i$ to $V_{g,n}(\LL)$. Therefore, one might wonder whether there are similar expressions for higher derivatives. In fact, we have the following equation involving the second derivative, although it is in some sense equivalent to the string equation.
\[
\frac{\partial^2 V_{g,n+1}}{\partial L_{n+1}^2}(\LL, 2\pi i) = \sum_{k=1}^n L_k \,\frac{\partial V_{g,n}(\LL)}{\partial L_k} - (4g-4+2n) \,V_{g,n}(\LL)
\]
There is reason to believe that such equations for higher derivatives simply do not exist. For example, see the work of Eynard and Orantin~\cite{eyn-ora}, which considers Weil--Petersson volumes as analogous to correlation functions arising from matrix models. They predict string and dilaton equations for functions which emerge from a vast generalisation of Mirzakhani's recursion.

\subsection{The large $g$ limit} \label{large-g}

Mirzakhani's recursion can in theory be used to calculate all Weil--Petersson volumes. However, a direct implementation of the recursion yields a computer program which is practical only for small genus. Zograf has provided an alternative algorithm which is empirically much faster~\cite{zog4}. In particular, he has managed to gather enough numerical evidence to suggest two interesting conjectures involving Weil--Petersson volumes in the large $g$ limit. In order to state the first, we use the following notation introduced by Mirzakhani~\cite{mir4} to express a certain normalisation of the coefficients of the polynomial $V_{g,n}(\LL)$.
\[
[\tau_{\alpha_1} \tau_{\alpha_2} \cdots \tau_{\alpha_n}]_{g,n} = \frac{\prod 2^{2\alpha_k} (2\alpha_k+1)!!}{(3g-3+n-|\bm{\alpha}|)!} \int_{\M_{g,n}} \psi_1^{\alpha_1} \psi_2^{\alpha_2} \cdots \psi_n^{\alpha_n} \omega^{3g-3+n-|\bm{\alpha}|}
\]

\begin{proposition} \label{zograf-2}
For a fixed tuple $(\alpha_1, \alpha_2, \ldots, \alpha_n)$ of non-negative integers,
\[
\lim_{g \to \infty} \frac{[\tau_{\alpha_1} \tau_{\alpha_2} \cdots \tau_{\alpha_n}]_{g,n}}{V_{g,n}(\mathbf{0})} = 1.
\]
\end{proposition}

Although stated as a conjecture by Zograf~\cite{zog4}, the result follows from certain Weil--Petersson volume estimates due to Mirzakhani~\cite{mir4}. These are obtained by rewriting her recursion --- see Theorem~\ref{mirzakhani-recursion} --- in the following form.
\begin{align*}
[\tau_{\alpha_1} \tau_{\alpha_2} \cdots \tau_{\alpha_n}]_{g,n} = & \frac{1}{2} \sum_{m=0}^{3g-3+n - |\bm{\alpha}|} \sum_{i + j = \alpha_1+m-2} b_m \,[\tau_i \tau_j \tau_{\alpha_2} \tau_{\alpha_3} \cdots \tau_{\alpha_n}]_{g-1,n+1} \\
& + \frac{1}{2} \sum_{\substack{I \sqcup J = [2,n] \\ g_1 + g_2 = g}} \sum_{m=0}^{3g-3+n - |\bm{\alpha}|} \sum_{i + j = \alpha_1+m-2} b_m \,[\tau_i \tau_{\bm{\alpha}_I}]_{g_1, |I|+1} \,[\tau_j \tau_{\bm{\alpha}_J}]_{g_2, |J|+1} \\
& + \sum_{k=2}^n \sum_{m=0}^{3g-3+n - |\bm{\alpha}|} (2\alpha_k+1) \,b_m \,[\tau_{\alpha_2} \cdots \tau_{\alpha_k+\alpha_1+m-1} \cdots \tau_{\alpha_n}]_{g,n-1}
\end{align*}
Here, we have used the notation $\tau_{\bm{\alpha}_I} = \tau_{\alpha_{i_1}} \tau_{\alpha_{i_2}} \cdots \tau_{\alpha_{i_m}}$ for $I = \{i_1, i_2, \ldots, i_m\}$. It is particularly useful to observe that the sequence $b_n = \zeta(2n)(1-\frac{1}{2^{2n-1}})$ consists only of positive terms, is strictly increasing, and limits to the value of 1.

The second of Zograf's conjectures gives the asymptotic behaviour of $V_{g,n}(\mathbf{0})$.

\begin{conjecture} \label{zograf-1}
For a fixed non-negative integer $n$,
\[
V_{g,n}(\mathbf{0}) = \frac{1}{\sqrt{g\pi}} (4\pi^2)^{2g-3+n} (2g-3+n)! \left[ 1 + c_ng^{-1} + O(g^{-2}) \right] \quad \text{as} \quad g \to \infty.
\]
\end{conjecture}

Although this conjecture is yet to be proven, Mirzakhani~\cite{mir4} has offered some supporting evidence.

\begin{theorem}
For a fixed non-negative integer $n$,
\[
\frac{V_{g,n+1}(\mathbf{0})}{2gV_{g,n}(\mathbf{0})} = 4\pi^2 + O(g^{-1}) \quad \text{and} \quad \frac{V_{g,n}(\mathbf{0})}{V_{g-1,n+2}(\mathbf{0})} = 1 + O(g^{-1}) \quad \text{as} \quad g \to \infty.
\]
\end{theorem}

Mirzakhani has used these and other Weil--Petersson volume estimates in the large $g$ limit to investigate the geometric properties of random hyperbolic surfaces~\cite{mir4}. In particular, she has obtained results concerning the length of the shortest simple closed geodesic, the diameter, and the Cheeger constant of a random surface with large genus, chosen with respect to the Weil--Petersson measure.

It is worth remarking that the $g \to \infty$ limit for fixed $n$ is much more difficult than the $n \to \infty$ limit for fixed $g$. The latter case was investigated by Manin and Zograf, who prove the following result~\cite{man-zog}.

\begin{theorem}
There exist constants $a_0, a_1, a_2, \ldots$ and $C$ such that, for a fixed non-negative integer $g$,
\[
V_{g,n}(\mathbf{0}) = n! C^n n^{(5g-7)/2} \left[ a_g + O(n^{-1}) \right] \quad \text{as} \quad n \to \infty.
\]
\end{theorem}

\subsection{The asymptotic Weil--Petersson form}

If we are only interested in psi-class intersection numbers on $\M_{g,n}$, then we need only look at the top degree part of the polynomial $V_{g,n}(\LL)$. This observation leads us to consider the following asymptotics of the Weil--Petersson volume for a fixed value of $\x = (x_1, x_2, \ldots, x_n)$.
\setcounter{equation}{0}
\begin{equation} \label{equation-1}
\lim_{N \to \infty} \frac{V_{g,n}(N\x)}{N^{6g-6+2n}} = \sum_{|\bm{\alpha}| = 3g-3+n} \frac{\int_{\M_{g,n}} \psi_1^{\alpha_1} \psi_2^{\alpha_2} \cdots \psi_n^{\alpha_n}}{2^{3g-3+n} \alpha_1! \alpha_2! \cdots \alpha_n!} \, x_1^{2\alpha_1} x_2^{2\alpha_2} \cdots x_n^{2\alpha_n}
\end{equation}

One way to access the asymptotics of the Weil--Petersson volume is via the following map on the combinatorial moduli space.
\[
f: \mathcal{MRG}_{g,n}(\x) \to \mathcal{MRG}_{g,n}(N\x) \to {\mathcal M}_{g,n}(N\x)
\]
This homeomorphism of orbifolds is the composition of two maps --- the first scales the ribbon graph metric by $N$ while the second uses the Bowditch--Epstein construction described in Section~\ref{combinatorial-moduli-space}. In one direction, this construction associates to a hyperbolic surface with boundary its spine --- in other words, the set of points which have at least two equal shortest paths to the boundary. The inverse of this construction produces a hyperbolic surface $S(\widetilde{\Gamma}) \in {\mathcal M}_{g,n}(\LL)$ for every metric ribbon graph $\widetilde{\Gamma} \in {\mathcal MRG}_{g,n}(\LL)$.

The normalised Weil--Petersson form $\frac{\omega}{N^2}$ on ${\mathcal M}_{g,n}(N\x)$ pulls back via $f$ to a symplectic form on the combinatorial moduli space. We will be interested in the limiting behaviour of this symplectic form since we may alternatively express the asymptotics of the Weil--Peterson volume in the following way.
\begin{align*}
\lim_{N \to \infty} \frac{V_{g,n}(N\x)}{N^{6g-6+2n}} &= \frac{1}{(3g-3+n)!}\lim_{N \to \infty} \int_{{\mathcal M}_{g,n}(N\x)} \left( \frac{\omega}{N^2} \right)^{3g-3+n} \\
&= \frac{1}{(3g-3+n)!} \int_{{\mathcal MRG}_{g,n}(\x)} \left( \lim_{N \to \infty} \frac{f^*\omega}{N^2} \right)^{3g-3+n} \\ 
&= \frac{1}{(3g-3+n)!} \sum_{\Gamma} \int_{{\mathcal MRG}_{\Gamma}(\x)} \left( \lim_{N \to \infty} \frac{f^*\omega}{N^2} \right)^{3g-3+n}
\end{align*}
To obtain the second line from the first, we pull back the integral to the combinatorial moduli space and invoke the Lebesgue dominated convergence theorem to move the limit inside the integral. To obtain the third line from the second, we use the orbifold cell decomposition of the combinatorial moduli space described in Section~\ref{combinatorial-moduli-space}. Recall that the combinatorial moduli space ${\mathcal MRG}_{g,n}(\x)$ is the disjoint union of open orbifold cells ${\mathcal MRG}_\Gamma(\x)$, where $\Gamma$ ranges over the ribbon graphs of type $(g,n)$. Here, the sum is only over the set of trivalent ribbon graphs of type $(g,n)$, since these correspond precisely to the open cells of this decomposition.

The previous discussion suggests that we should study the asymptotic behaviour of the Weil--Petersson form. In order to do this, fix a trivalent ribbon graph $\Gamma$ of type $(g,n)$ and label its edges from 1 up to $6g-6+3n$. As noted in Section~\ref{combinatorial-moduli-space}, the lengths of these edges $e_1, e_2, \ldots, e_{6g-6+3n}$ provide a set of natural coordinates on ${\mathcal MRG}_\Gamma(\x)$ and we can write
\[
{\mathcal MRG}_\Gamma(\x) \cong \left. \left\{ \mathbf{e} \in \mathbb{R}_+^{6g-6+3n} \; \middle \vert \; A_\Gamma \mathbf{e} = \x \right\} \right/ \mathrm{Aut}(\Gamma).
\]
Here, $A_\Gamma$ is the linear map which represents the adjacency between faces and edges in the cell decomposition corresponding to $\Gamma$.

\begin{theorem} \label{equal-forms}
In the $N \to \infty$ limit, the symplectic form $\frac{f^*\omega}{N^2}$ converges pointwise to a 2-form $\Omega$ on ${\mathcal MRG}_\Gamma(\x)$. There exists a $(6g-6+2n) \times (6g-6+2n)$ skew-symmetric matrix $B_\Gamma$ such that, after an appropriate permutation of the edge labels,
\[
\Omega = \sum_{1 \leq i < j \leq 6g-6+2n} (B_\Gamma)_{ij} \, de_i \wedge de_j.
\]
\end{theorem}

We can use this theorem to write the asymptotics of the Weil--Petersson volume in the following way, where $\mathrm{pf}(B_\Gamma)$ denotes the Pfaffian of $B_\Gamma$.
\setcounter{equation}{1}
\begin{equation} \label{equation-2}
\lim_{N \to \infty} \frac{V_{g,n}(N\x)}{N^{6g-6+2n}} = \sum_{\Gamma} \frac{\mathrm{pf}(B_\Gamma)}{|\mathrm{Aut}(\Gamma)|} \int_{A_\Gamma \mathbf{e} = \mathbf{x}} de_1 \wedge de_2 \wedge \cdots \wedge de_{6g-6+2n}
\end{equation}
Now we can equate the expressions appearing in Equations~\eqref{equation-1} and \eqref{equation-2}. The Pfaffian and integral in Equation~\eqref{equation-2} can be calculated explicitly in terms of the combinatorics of $\Gamma$. Upon doing so and taking the Laplace transform of both sides, we recover the following identity.

\begin{theorem}[Kontsevich's combinatorial formula] \label{kcf} For non-negative $g$ and positive $n$ satisfying $2 - 2g - n < 0$, we have the following equality of rational polynomials in $s_1, s_2, \ldots, s_n$.
\[
\sum_{|\bm{\alpha}| = 3g-3+n} \int_{\M_{g,n}} \psi_1^{\alpha_1} \psi_2^{\alpha_2} \cdots \psi_n^{\alpha_n} \prod_{k=1}^n \frac{(2\alpha_k - 1)!!}{s_k^{2\alpha_k+1}} = \sum_{\Gamma} \frac{2^{2g-2+n}}{|\mathrm{Aut}(\Gamma)|} \prod_{e \in E(\Gamma)} \frac{1}{s_{\ell(e)} + s_{r(e)}}
\]
The sum on the right hand side is over trivalent ribbon graphs of type $(g,n)$. For an edge $e$, the expressions $\ell(e)$ and $r(e)$ denote the labels of the faces on its left and right.\footnote{Although the left and right of an edge are not well-defined, the expression $s_{\ell(e)} + s_{r(e)}$ certainly is.}
\end{theorem}

This is the main identity used by Kontsevich in his proof of the Witten--Kontsevich theorem. Our proof of this result highlights the close relationship between the combinatorial methods pioneered by Kontsevich \cite{kon} and the hyperbolic geometry used by Mirzakhani \cite{mir2, mir1}.

It is worth making a few remarks on the proof of Theorem~\ref{equal-forms}. The result appears implicitly in the work of Mondello~\cite{mon} although we will discuss an alternative proof which is less computational in nature~\cite{do2}. Underlying this work is the observation that a hyperbolic surface with large boundary lengths resembles a ribbon graph after appropriately scaling the hyperbolic metric. In order to make this statement precise, take a metric ribbon graph $\widetilde{\Gamma}$ and consider the surface $\frac{1}{N} S(N\widetilde{\Gamma})$ for large values of $N$. Here, we use the notation $\lambda X$ to denote the result of scaling the metric on $X$ by a positive real number $\lambda$. By the Gauss--Bonnet theorem, the area of the surface $\frac{1}{N} S(N\widetilde{\Gamma})$ goes to zero as $N$ increases to infinity. On the other hand, the length of the boundaries remains fixed. So in the $N \to \infty$ limit, one expects the entire surface to collapse onto the spine $\widetilde{\Gamma}$. The following result formalises this intuition in a precise way~\cite{do2}.

\begin{theorem} \label{gromov-hausdorff}
In the Gromov--Hausdorff topology, for every metric ribbon graph $\widetilde{\Gamma}$, we have
\[
\lim_{N \to \infty} \frac{1}{N} S(N\widetilde{\Gamma})= \widetilde{\Gamma}.
\]
\end{theorem}

The intuitive observation behind this result suggests a great deal about the geometry of a hyperbolic surface with large boundary lengths. For example, one expects the length of a closed geodesic on $\frac{1}{N} S(N\widetilde{\Gamma})$  to converge to a sum of lengths of edges in $\widetilde{\Gamma}$. Furthermore, one expects the acute angle at which two closed geodesics meet to converge to 0. We recover the limiting behaviour of the Weil--Petersson form from the limiting behaviour of such lengths and angles via the following result of Wolpert~\cite{wol1}.

\begin{proposition} \label{wolpert}
Let $C_1, C_2, \ldots, C_{6g-6+2n}$ be simple closed geodesics with lengths $\ell_1, \ell_2, \ldots, \ell_{6g-6+2n}$ in a hyperbolic surface $S \in {\mathcal M}_{g,n}(\LL)$. If $C_i$ and $C_j$ meet at a point $p$, denote by $\theta_p$ the angle between the curves, measured anticlockwise from $C_i$ to $C_j$. Define the $(6g-6+2n) \times (6g-6+2n)$ skew-symmetric matrix $X$ by the formula
\[
X_{ij} = \sum_{p \in C_i \cap C_j} \cos \theta_p, \qquad \text{for } i < j.
\]
If $X$ is invertible, then $\ell_1, \ell_2, \ldots, \ell_{6g-6+2n}$ are local coordinates at $S \in {\mathcal M}_{g,n}(\LL)$ and the Weil--Petersson form is given by
\[
\omega = - \sum_{i < j}~[X^{-1}]_{ij}~d\ell_i \wedge d\ell_j.
\]
\end{proposition}

A judicious choice of curves allows us to use this result to obtain Theorem~\ref{equal-forms}, including a concrete description of the matrix $B_\Gamma$ in terms of the combinatorics of the ribbon graph $\Gamma$. Furthermore, one finds that the asymptotic Weil--Petersson form $\Omega$ coincides with the 2-form on the combinatorial moduli space introduced by Kontsevich in his proof of the Witten--Kontsevich theorem~\cite{kon}. For the details of the proof of Theorem~\ref{equal-forms}, consult the relevant source in the literature~\cite{do2}.

\appendix

\section{Intersection theory on moduli spaces of curves} \label{intersection-theory}

The reader will find a wealth of information concerning moduli spaces of curves and their intersection theory elsewhere in the literature~\cite{har-mor, vak1, vak2}. The aim of this appendix is to provide a concise exposition of the topic in order to keep this article reasonably self-contained.

For non-negative integers $g$ and $n$ satisfying the Euler characteristic condition $2 - 2g - n < 0$, define the {\em moduli space of curves} as follows.
\[
{\mathcal M}_{g,n} = \left. \left\{ (C, p_1, p_2, \ldots, p_n) \; \middle \vert \; \begin{array}{l} C \text{ is a smooth algebraic curve with genus } g \\ \text{and } n \text{ distinct points } p_1, p_2, \ldots, p_n \end{array} \right\} \right / \sim
\]
Here, $(C, p_1, p_2, \ldots, p_n) \sim (D, q_1, q_2, \ldots, q_n)$ if and only if there exists an isomorphism from $C$ to $D$ which sends $p_k$ to $q_k$ for all $k$. It is often more natural to work with the {\em Deligne--Mumford compactification} of the moduli space of curves.
\[
\overline{\mathcal M}_{g,n} = \left. \left\{ (C, p_1, p_2, \ldots, p_n) \; \middle \vert \; \begin{array}{l} C \text{ is a stable algebraic curve with genus } g \\ \text{and } n \text{ distinct smooth points } p_1, p_2, \ldots, p_n \end{array} \right\} \right / \sim
\]
Again, $(C, p_1, p_2, \ldots, p_n) \sim (D, q_1, q_2, \ldots, q_n)$ if and only if there exists an isomorphism from $C$ to $D$ which sends $p_k$ to $q_k$ for all $k$. An algebraic curve is called {\em stable} if it has at worst nodal singularities and a finite automorphism group. The practical interpretation of this latter condition is that the normalisation of every rational component must have at least three distinguished points which are nodes or marked points. One of the virtues of the Deligne--Mumford compactification amongst the various competing options is the fact that it is modular --- in other words, each point in $\M_{g,n}$ represents an algebraic curve. The set $\M_{g,n}$ possesses a rich geometric structure and is an example of a Deligne--Mumford stack, although one will not go too far wrong thinking of it as a complex orbifold.

A natural approach to understanding the structure of geometric spaces is through algebraic invariants, such as homology and cohomology. And so it is with moduli spaces of curves, but for the fact that its full cohomology ring is notoriously intractable in general. However, a great deal of progress can be made by calculating intersection numbers with respect to certain characteristic classes. The classes that we consider live in the cohomology ring $H^*(\M_{g,n}; \mathbb{Q})$ and arise from taking Chern classes of natural complex vector bundles.\footnote{Readers with a more algebraic predilection may prefer to think of these classes as living in the Chow ring $A^*(\M_{g,n})$.} One obtains a more natural theory using rational, rather than integral, coefficients for cohomology due to the orbifold nature of $\M_{g,n}$.

Given a stable genus $g$ curve with $n+1$ labelled points, one can forget the point labelled $n+1$ to obtain a genus $g$ curve with $n$ labelled points. The resulting curve may not be stable, but gives rise to a well-defined stable curve after contracting all unstable rational components. This yields a map $\pi: \M_{g,n+1} \to \M_{g,n}$ known as the {\em forgetful morphism}, which can be interpreted as the universal family over $\M_{g,n}$. Thus, given a pair $(C,p)$ consisting of a stable curve $C \in \M_{g,n}$ and a point $p$ on the curve, it is possible to associate to it a unique stable curve $D \in \M_{g,n+1}$ such that $\pi(D) = C$. In particular, the fibre over $C \in \M_{g,n}$ is essentially the stable curve corresponding to $C$. So the point labelled $k$ defines a section $\sigma_k: \M_{g,n} \to \M_{g,n+1}$ for $k = 1, 2, \ldots, n$. The forgetful morphism can be used to pull back cohomology classes, but it will also be useful to push them forward. This is possible via the Gysin map $\pi_*: H^*(\M_{g,n+1}; \mathbb{Q}) \to H^*(\M_{g,n}; \mathbb{Q})$, the homomorphism of graded rings with grading $-2$ which represents integration along fibres.

Consider the vertical cotangent bundle on $\M_{g,n+1}$ whose fibre at the point associated to the pair $(C,p)$ is equal to the cotangent line $T_p^*C$. Unfortunately, this definition is nonsensical when $p$ is a singular point of $C$. Therefore, it is necessary to consider the relative dualising sheaf, the unique line bundle on $\M_{g,n+1}$ which extends the vertical cotangent bundle. More precisely, it can be defined as ${\mathcal L} = {\mathcal K}_X \otimes \pi^* {\mathcal K}_B^{-1}$, where ${\mathcal K}_X$ denotes the canonical line bundle on $\M_{g,n+1}$ and ${\mathcal K}_B$ denotes the canonical line bundle on $\M_{g,n}$. Sections of ${\mathcal L}$ correspond to meromorphic 1-forms with at worst simple poles allowed at the nodes which also satisfy the condition that the two residues at the preimages of each node under normalisation must sum to zero.

The tautological line bundles on $\M_{g,n}$ are formed by pulling back $\mathcal L$ along the sections $\sigma_k$ for $k = 1, 2, \ldots, n$. Taking Chern classes of these line bundles, we obtain the {\em psi-classes}
\[
\psi_k = c_1(\sigma_k^* {\mathcal L}) \in H^2(\M_{g,n}; \mathbb{Q}) \qquad \text{for } k = 1, 2, \ldots, n.
\]

Define the twisted Euler class by $e = c_1 \left({\mathcal L} \left( D_1 + D_2 + \cdots + D_n \right) \right)$, where $D_k$ is the divisor on $\M_{g,n+1}$ representing the image of the section $\sigma_k$. Taking the pushforwards of its powers, we obtain the {\em Mumford--Morita--Miller classes}
\[
\kappa_m = \pi_*(e^{m+1}) \in H^{2m}(\M_{g,n}; \mathbb{Q}) \qquad \text{for } m = 0, 1, 2, \ldots, 3g-3+n.
\]

A great deal of attention has been paid to the subring of $H^*(\M_{g,n}; \mathbb{Q})$ known as the {\em tautological ring}. It has the benefit of being more tractable than the full cohomology ring and possessing a rich combinatorial structure, while still containing all known classes of geometric interest. Any top intersections in the tautological ring can be determined from the top intersections of psi-classes alone. Thus, we are motivated to study intersection numbers of the form
\[
\langle \tau_{\alpha_1} \tau_{\alpha_2} \cdots \tau_{\alpha_n} \rangle = \int_{\M_{g,n}} \psi_1^{\alpha_1} \psi_2^{\alpha_2} \cdots \psi_n^{\alpha_n} \in \mathbb{Q},
\]
where $|\bm{\alpha}| = 3g - 3 + n$ or equivalently, $g = \frac{1}{3}(|\bm{\alpha}| - n + 3)$. The bracket notation above --- originally introduced by Witten --- suppresses the genus and encodes the symmetry between the psi-classes. We treat the $\tau$ variables as commuting, so that we can write intersection numbers in the form $\langle \tau_0^{d_0} \tau_1^{d_1} \tau_2^{d_2} \cdots \rangle$ and we set $\langle \tau_{\alpha_1} \tau_{\alpha_2} \cdots \tau_{\alpha_n} \rangle = 0$ if $n = 0$ or if the genus $g = \frac{1}{3}(|\bm{\alpha}| - n + 3)$ is non-integral or negative. In this way, we have defined a linear functional $\langle \cdot \rangle: \mathbb{Q}[\tau_0, \tau_1, \tau_2, \ldots] \to \mathbb{Q}$. The psi-class intersection numbers contain a great deal of structure, as evidenced by the following result~\cite{wit}.

\begin{proposition}[String and dilaton equations] \label{string-dilaton} For $2g-2+n > 0$, the psi-class intersection numbers satisfy the following relations.
\begin{align*}
\langle \tau_0 \tau_{\alpha_1} \tau_{\alpha_2} \cdots \tau_{\alpha_n} \rangle &= \sum_{k=1}^n \langle \tau_{\alpha_1} \cdots \tau_{\alpha_k-1} \cdots \tau_{\alpha_n} \rangle \\
\langle \tau_1 \tau_{\alpha_1} \tau_{\alpha_2} \cdots \tau_{\alpha_n} \rangle &= (2g-2+n) \langle
\tau_{\alpha_1} \tau_{\alpha_2} \cdots \tau_{\alpha_n} \rangle
\end{align*}
\end{proposition}

One of the landmark results concerning intersection theory on moduli spaces of curves is Witten's conjecture, now Kontsevich's theorem. In conjunction with the string equation, it allows us to calculate any psi-class intersection number from the base case $\langle \tau_0^3 \rangle = 1$. In order to precisely describe the result, we let $\mathbf{t} = (t_0, t_1, t_2, \ldots)$ and $\mathbf{\tau} = (\tau_0, \tau_1, \tau_2, \ldots)$ and consider the generating function $F(\mathbf{t}) = \langle \exp(\mathbf{t} \cdot \mathbf{\tau}) \rangle$. Here, the expression is to be expanded using multilinearity in the variables $t_0, t_1, t_2, \ldots$. Equivalently, we may define
\[
F(t_0, t_1, t_2, \ldots) = \sum_{\mathbf{d}} \prod_{k=0}^\infty \frac{t_k^{d_k}}{d_k!} \langle \tau_0^{d_0} \tau_1^{d_1} \tau_2^{d_2} \cdots \rangle,
\]
where the summation is over all sequences $\mathbf{d} = (d_0, d_1, d_2, \ldots)$ of non-negative integers with finitely many non-zero terms. In his foundational paper~\cite{wit}, Witten argued on physical grounds that the formal series $U = \frac{\partial^2 F}{\partial t_0^2}$ satisfies the KdV hierarchy of partial differential equations. This is the prototypical example of an exactly solvable model, whose soliton solutions have attracted tremendous mathematical interest over the past few decades. More explicitly, the Witten--Kontsevich theorem can be stated in the following way.

\begin{theorem}[Witten--Kontsevich theorem]
The generating function $F$ satisfies the following partial differential equation for every non-negative integer $n$.
\[
(2n+1) \frac{\partial^3 F}{\partial t_n \partial t_0^2} = \left( \frac{\partial^2 F}{\partial t_{n-1} \partial t_0} \right) \left( \frac{\partial^3 F}{\partial t_0^3} \right) + 2 \left( \frac{\partial^3 F}{\partial t_{n-1} \partial t_0^2} \right) \left( \frac{\partial^2 F}{\partial t_0^2} \right) + \frac{1}{4} \frac{\partial^5 F}{\partial t_{n-1} \partial t_0^4}
\]
\end{theorem}

An equivalent formulation of the Witten--Kontsevich theorem states that the Virasoro operators annihilate the generating function $\exp F$. These operators span the Virasoro Lie algebra and are defined by
\[
L_{-1} = - \frac{1}{2} \frac{\partial}{\partial t_0} + \frac{1}{2} \sum_{k=0}^\infty t_{k+1} \frac{\partial}{\partial t_k} + \frac{t_0^2}{4}, \quad L_0 = - \frac{3}{2} \frac{\partial}{\partial t_1} + \frac{1}{2} \sum_{k=0}^\infty (2k+1) t_k \frac{\partial}{\partial t_k} + \frac{1}{48},
\]
and for positive integers $n$,
\[
L_n = -\frac{(2n+3)!!}{2} \frac{\partial}{\partial t_{n+1}} + \sum_{k=0}^\infty \frac{(2k+2n+1)!!}{2(2k-1)!!} t_k \frac{\partial}{\partial t_{n+k}} + \sum_{i + j = n-1} \frac{(2i+1)!! (2j+1)!!}{4} \frac{\partial^2}{\partial t_i \partial t_j}.
\]

There now exist several proofs of the Witten--Kontsevich theorem, due to Kontsevich~\cite{kon}, Okounkov and Pandharipande~\cite{oko-pan}, Kim and Liu~\cite{kim-liu}, Kazarian and Lando~\cite{kaz-lan}, and Mirzakhani~\cite{mir1}. That there are so many proofs, each with their own distinct flavour, is testament to the importance and richness of the result.

\section{Table of Weil--Petersson volumes} \label{WP-table}

The following table shows some examples of Weil--Petersson volumes. We use the notation $m_{(\alpha_1, \alpha_2, \ldots, \alpha_k)}$ to denote the monomial symmetric polynomial
\[
\sum_{(\beta_1, \beta_2, \ldots, \beta_n)} L_1^{2\beta_1} L_2^{2\beta_2} \cdots L_n^{2\beta_n},
\]
where the summation ranges over all permutations $(\beta_1, \beta_2, \ldots, \beta_n)$ of $(\alpha_1, \alpha_2, \ldots, \alpha_k, 0, 0, \ldots, 0)$. For example, we have the following when $n = 3$.
\begin{align*}
m_{(3,2,1)} &= L_1^6L_2^4L_3^2 + L_1^6L_2^2L_3^4 + L_1^4L_2^6L_3^2 + L_1^4L_2^2L_3^6 + L_1^2L_2^6L_3^4 + L_1^2L_2^4L_3^6 \\ m_{(2)} &= L_1^4 + L_2^4 + L_3^4 \\
m_{(1,1,1)} &= L_1^2 L_2^2 L_3^2
\end{align*}

\begin{tabular*}{\textwidth}{@{}ccl@{}} \toprule
$g$ & $n$ & $V_{g,n}(L_1, L_2, \ldots, L_n)$ \\ \midrule
0 & 3 & 1 \\
 & 4 & $\frac{1}{2} m_{(1)} + 2\pi^2$ \\
 & 5 & $\frac{1}{8} m_{(2)} + \frac{1}{2} m_{(1,1)} + 3\pi^2 m_{(1)} + 10\pi^4$ \\
 & 6 & $\frac{1}{48} m_{(3)} + \frac{3}{16} m_{(2,1)} + \frac{3}{4} m_{(1,1,1)} + \frac{3\pi^2}{2} m_{(2)} + 6\pi^2 m_{(1,1)} + 26\pi^4 m_{(1)} + \frac{244\pi^6}{3}$ \\ 
 & 7 & $\frac{1}{384} m_{(4)} + \frac{1}{24} m_{(3,1)} + \frac{3}{32} m_{(2,2)} + \frac{3}{8} m_{(2,1,1)} + \frac{3}{2} m_{(1,1,1,1)} + \frac{5\pi^2}{12} m_{(3)} + \frac{15\pi^2}{4} m_{(2,1)} + 15\pi^2 m_{(1,1,1)}$ \\
 & & $+ 20\pi^4 m_{(2)} + 80\pi^4 m_{(1,1)} + \frac{910\pi^6}{3} m_{(1)} + \frac{2758 \pi^8}{3}$ \\ \bottomrule
 \end{tabular*}

\begin{tabular*}{\textwidth}{@{}ccl@{}} \toprule
$g$ & $n$ & $V_{g,n}(L_1, L_2, \ldots, L_n)$ \\ \midrule
1 & 1 & $\frac{1}{48} m_{(1)} + \frac{\pi^2}{12}$ \\
 & 2 & $\frac{1}{192} m_{(2)} + \frac{1}{96} m_{(1,1)} + \frac{\pi^2}{12} m_{(1)} + \frac{\pi^4}{4}$ \\
 & 3 & $\frac{1}{1152} m_{(3)} + \frac{1}{192} m_{(2,1)} + \frac{1}{96} m_{(1,1,1)} + \frac{\pi^2}{24} m_{(2)} + \frac{\pi^2}{8} m_{(1,1)} + \frac{13\pi^4}{24} m_{(1)} + \frac{14\pi^6}{9}$ \\
 & 4 & $\frac{1}{9216} m_{(4)} + \frac{1}{768} m_{(3,1)} + \frac{1}{384} m_{(2,2)} + \frac{1}{128} m_{(2,1,1)} + \frac{1}{64} m_{(1,1,1,1)} + \frac{7\pi^2}{576} m_{(3)} + \frac{\pi^2}{12} m_{(2,1)}$ \\
 & & $ + \frac{\pi^2}{4} m_{(1,1,1)} + \frac{41\pi^4}{96} m_{(2)} + \frac{17\pi^4}{12} m_{(1,1)} + \frac{187\pi^6}{36} m_{(1)} + \frac{529\pi^8}{36}$ \\
 & 5 & $\frac{1}{92160} m_{(5)} + \frac{1}{4608} m_{(4,1)} + \frac{7}{9216} m_{(3,2)} + \frac{1}{384} m_{(3,1,1)} + \frac{1}{192} m_{(2,2,1)} + \frac{1}{64} m_{(2,1,1,1)} + \frac{1}{32} m_{(1,1,1,1,1)}$ \\
 & & $+ \frac{11\pi^2}{4608} m_{(4)} + \frac{35\pi^2}{1152} m_{(3,1)} + \frac{\pi^2}{16} m_{(2,2)} + \frac{5\pi^2}{24} m_{(2,1,1)} + \frac{5\pi^2}{8} m_{(1,1,1,1)} + \frac{13\pi^4}{72} m_{(3)} + \frac{253\pi^4}{192} m_{(2,1)}$ \\
 & & $+ \frac{35\pi^4}{8} m_{(1,1,1)} + \frac{809\pi^6}{144} m_{(2)} + \frac{703\pi^6}{36} m_{(1,1)} + \frac{4771\pi^8}{72} m_{(1)} + \frac{16751\pi^{10}}{90}$ \\ \midrule
2 & 0 & $\frac{43\pi^6}{2160}$ \\
 & 1 & $\frac{1}{442368} m_{(4)} + \frac{29\pi^2}{138240} m_{(3)} + \frac{139\pi^4}{23040} m_{(2)} + \frac{169\pi^6}{2880} m_{(1)} + \frac{29\pi^8}{192}$ \\
 & 2 & $\frac{1}{4423680} m_{(5)} + \frac{1}{294912} m_{(4,1)} + \frac{29}{2211840} m_{(3,2)} + \frac{11\pi^2}{276480} m_{(4)} + \frac{29\pi^2}{69120} m_{(3,1)} + \frac{7\pi^2}{7680} m_{(2,2)} + \frac{19\pi^4}{7680} m_{(3)}$ \\
 & & $+ \frac{181\pi^4}{11520} m_{(2,1)} + \frac{551\pi^6}{8640} m_{(2)} + \frac{7\pi^6}{36} m_{(1,1)} + \frac{1085\pi^8}{1728} m_{(1)} + \frac{787\pi^{10}}{480}$ \\
& 3 & $\frac{1}{53084160} m_{(6)} + \frac{1}{2211840} m_{(5,1)} + \frac{11}{4423680} m_{(4,2)} + \frac{1}{147456} m_{(4,1,1)} + \frac{29}{6635520} m_{(3,3)} + \frac{29}{1105920} m_{(3,2,1)}$ \\
 & & $+ \frac{7}{122880} m_{(2,2,2)} + \frac{\pi^2}{172800} m_{(5)} + \frac{11\pi^2}{110592} m_{(4,1)} + \frac{5\pi^2}{13824} m_{(3,2)} + \frac{29\pi^2}{27648} m_{(3,1,1)} + \frac{7\pi^2}{3072} m_{(2,2,1)}$ \\
 & & $+ \frac{41\pi^4}{61440} m_{(4)} + \frac{211\pi^4}{27648} m_{(3,1)} + \frac{37\pi^4}{2304} m_{(2,2)} + \frac{223\pi^4}{4608} m_{(2,1,1)} + \frac{77\pi^6}{2160} m_{(3)} + \frac{827\pi^6}{3456} m_{(2,1)}$ \\
 & & $+ \frac{419\pi^6}{576} m_{(1,1,1)} + \frac{30403\pi^8}{34560} m_{(2)} + \frac{611\pi^8}{216} m_{(1,1)} + \frac{75767\pi^{10}}{8640} m_{(1)} + \frac{1498069\pi^{12}}{64800}$ \\ \midrule
3 & 0 & $\frac{176557\pi^{12}}{1209600}$ \\
 & 1 & $\frac{1}{53508833280} m_{(7)} + \frac{77\pi^2}{9555148800} m_{(6)} + \frac{3781\pi^4}{2786918400} m_{(5)} + \frac{47209\pi^6}{418037760} m_{(4)} + \frac{127189\pi^8}{26127360} m_{(3)}$ \\
 & & $+ \frac{8983379\pi^{10}}{87091200} m_{(2)} + \frac{8497697\pi^{12}}{9331200} m_{(1)} + \frac{9292841\pi^{14}}{4082400}$ \\ 
 & 2 & $\frac{1}{856141332480} m_{(8)} + \frac{1}{21403533312 } m_{(7,1)} + \frac{77}{152882380800} m_{(6,2)} + \frac{503}{267544166400} m_{(5,3)}$ \\
 & & $+ \frac{607}{214035333120} m_{(4,4)} + \frac{17\pi^2}{22295347200} m_{(7)} + \frac{77\pi^2}{3185049600} m_{(6,1)} + \frac{17\pi^2}{88473600} m_{(5,2)} + \frac{1121\pi^2}{2229534720} m_{(4,3)}$ \\
 & & $+ \frac{1499\pi^4}{7431782400} m_{(6)} + \frac{899\pi^4}{185794560} m_{(5,1)} + \frac{10009\pi^4}{371589120} m_{(4,2)} + \frac{191\pi^4}{4128768} m_{(3,3)} + \frac{3859\pi^6}{139345920} m_{(5)}$ \\
 & & $+ \frac{33053\pi^6}{69672960} m_{(4,1)} + \frac{120191\pi^6}{69672960} m_{(3,2)} + \frac{195697\pi^8}{92897280} m_{(4)} + \frac{110903\pi^8}{4644864} m_{(3,1)} + \frac{6977\pi^8}{138240} m_{(2,2)}$ \\
 & & $+ \frac{37817\pi^{10}}{430080} m_{(3)} + \frac{2428117\pi^{10}}{4147200} m_{(2,1)} + \frac{5803333\pi^{12}}{3110400} m_{(2)} + \frac{18444319\pi^{12}}{3110400} m_{(1,1)} + \frac{20444023\pi^{14}}{1209600} m_{(1)}$ \\
 & & $+ \frac{2800144027\pi^{16}}{65318400}$ \\ \midrule
4 & 0 & $\frac{1959225867017\pi^{18}}{493807104000}$ \\
 & 1 & $\frac{1}{29588244450508800} m_{(10)} + \frac{149\pi^2}{3698530556313600} m_{(9)} + \frac{48689\pi^4}{2397195730944000} m_{(8)} + \frac{50713\pi^6}{8989483991040} m_{(7)}$ \\
 & & $+ \frac{30279589\pi^8}{32105299968000} m_{(6)} + \frac{43440449\pi^{10}}{445906944000} m_{(5)} + \frac{274101371\pi^{12}}{44590694400} m_{(4)} + \frac{66210015481\pi^{14}}{292626432000} m_{(3)}$ \\
 & & $+ \frac{221508280867\pi^{16}}{50164531200} m_{(2)} + \frac{74706907467169\pi^{18}}{1975228416000} m_{(1)} + \frac{92480712720869\pi^{20}}{987614208000}$ \\ \midrule
5 & 0 & $\frac{84374265930915479\pi^{24}}{355541114880000}$ \\
 & 1 & $\frac{1}{48742490377990176768000} m_{(13)} + \frac{7\pi^2}{133907940598874112000} m_{(12)} + \frac{1823\pi^4}{31067656673034240000} m_{(11)}$ \\
 & & $+ \frac{296531\pi^6}{7766914168258560000} m_{(10)} + \frac{68114707\pi^8}{4271802792542208000} m_{(9)} + \frac{2123300941\pi^{10}}{474644754726912000} m_{(8)}$ \\
 & & $+ \frac{42408901133\pi^{12}}{49442161950720000} m_{(7)} + \frac{19817320001\pi^{14}}{176579149824000} m_{(6)} + \frac{11171220559409\pi^{16}}{1135151677440000} m_{(5)} + \frac{62028372646367\pi^{18}}{111244864389120} m_{(4)}$ \\
 & & $+ \frac{202087901261599\pi^{20}}{10534551552000} m_{(3)} + \frac{626693680890100121\pi^{22}}{1738201006080000} m_{(2)} + \frac{881728936440038779\pi^{24}}{289700167680000} m_{(1)}$ \\
 & & $+ \frac{21185241498983729441\pi^{26}}{2824576634880000}$ \\ \bottomrule
\end{tabular*}

\addcontentsline{toc}{section}{References}
\bibliographystyle{hacm}
\bibliography{Do-Handbook}

\end{document}